\documentclass[reqno]{amsart}
\newtheorem{thm}{Theorem}[section]
\newtheorem{df}[thm]{Definition}
\newtheorem{prop}[thm]{Proposition}
\newtheorem{lem}[thm]{Lemma}
\newtheorem{cor}[thm]{Corollary}
\newtheorem{exm}[thm]{Example}
\newtheorem{rem}[thm]{Remark}
\usepackage{amsmath,amssymb,amscd,epsfig}
\usepackage{color, graphics}
\numberwithin{equation}{section}

\newcommand{\wed}{\wedge}

\newcommand{\ag}{(\alpha,a)}

\newcommand{\GG}{\Gamma\times {\mathbb Z}_2 }

\newcommand{\ld}{\lambda}
\newcommand{\ch}{\text{ch}}

\newcommand{\pr}{\partial}

\newcommand{\bt}{\beta}

\newcommand{\G}{\Gamma}
\newcommand{\al}{\alpha}
\newcommand{\C}{\mathbb C}
\newcommand{\Z}{\mathbb Z}

\newcommand{\fL} {\frak L}

\newcommand{\la} {\lambda}

\begin{document}
\title[NECKLACE RINGS]
{NECKLACE RINGS AND LOGARITHMIC FUNCTIONS}
\author{YOUNG-TAK OH}
\address{Korea Institute for Advanced Study\\
         207-43 Cheongryangri-dong, Dongdaemun-gu\\
         Seoul 130-722, Korea}
\email{ohyt@kias.re.kr} \maketitle \baselineskip=12pt

\begin{abstract}
In this paper, we develop the theory of the necklace ring and the
logarithmic function. Regarding the necklace ring, we introduce
the necklace ring functor $Nr$ from the category of special
$\ld$-rings into the category of special $\ld$-rings and then
study the associated Adams operators. As far as the logarithmic
function is concerned, we generalize the results in Bryant's paper
(J. Algebra. 253 (2002); no.1, 167-188) to the case of graded Lie
(super)algebras with a group action by applying the
Euler-Poincar\'e principle.
\end{abstract}

\section{Introduction}
\renewcommand{\thefootnote}{}
\footnotetext[1]{%
\renewcommand{\baselinestretch}{1.2}\selectfont
This research was supported by KOSEF Grant \#
R01-2003-000-10012-0. \hfill \break MSC : 11F03,11F22,17B70.
\hfill \break
Keywords: Necklace ring, Logarithmic function, Lie superalgebra}

Let $V$ be a finite-dimensional vector space over an arbitrary field.
Many mathematicians have long intensively studied
the free Lie algebra $\mathfrak L(V)$ generated by $V$
for its remarkable connections with combinatorics.
In this paper, we study two algebraic objects
which originate from the combinatorics of free Lie algebras.
The necklace ring is related to the
well-known fact that the dimension of the $n$-th homogeneous
component in $\mathfrak L(V)$, can be computed by counting the
number of primitive necklaces of length $n$ out of $\dim V$ letters.
The second, the logarithmic function, is related to
the ``Lazard elimination theorem'' (\cite {Bour}).

The notion of the necklace ring was first introduced
by Metropolis and Rota.
The necklace ring has many interesting and significant
algebraic features.
For example, it is isomorphic to the universal ring of Witt vectors
over a certain class of commutative rings (\cite{MR}).
In particular, the necklace ring over $\mathbb Z$
was explicitly realized as
the Burnside-Grothendieck ring, $\hat \Omega(C)$,
of isomorphism classes of almost finite cyclic sets.
Here, the notation $C$ represents the multiplicative infinite cyclic group.
The isomorphism between the universal ring of Witt vectors,
$\mathbb W(\mathbb Z)$, and $\hat \Omega(C)$
has been called the {\it Teichm\"{u}ller map} (\cite{DS}).

The classical construction of Witt vectors can be understood
as a special case of a more general construction.
More precisely, Dress and Siebeneicher
constructed a covariant functor, $\mathbb W_G$,
such that $\mathbb W_{\hat C}$ coincides
with the classical Witt-ring functor, $\mathbb W$,
where $C$ denotes a profinite group and $\hat C$ the profinite
completion of $C$.
In this case, the isomorphism between
$\mathbb W_G(\mathbb Z)$ and $\hat \Omega(G)$ has been called the
{\it extended Teichm\"{u}ller map} (\cite{DS2}).
Generalizing this, Graham  constructed a functor, $F_G$,
which shares many properties with $\mathbb W_G$, for every group $G$.
Recently, Brun showed that for any finite group, $G$,
the functor $\mathbb W_G$ coincides with the left adjoint of the algebraic functor
from the category of $G$-Tambara functors to the category of
commutative rings with an action of $G$ (\cite{B,Gr}).

Motivated by the work of Dress and Siebeneicher,
we introduced a new ring, $\hat \Omega_R(G)$,
which coincides with $\hat \Omega(G)$, if $G=\hat C$
for a special $\lambda$-ring $R$.
Furthermore, we constructed a map,
$$\tau_R :\mathbb W_G(R) \to \hat \Omega_R(G),$$
which is analogous to the extended Teichm\"{u}ller map.
This map, which is a ring isomorphism in cases where $R$ is torsion-free,
has been called the {\it $R$-Teichm\"{u}ller map}
(\cite{O}).

The necklace ring of $G$ over $R$, denoted by $Nr_G(R)$, is the ring
obtained from $\hat \Omega_R(G)$ via the interpretation map (int for short)
which is nothing but the map reading of the coefficients
of the elements in $\hat \Omega_R(G)$.
Pictorially, this result can be illustrated in the
following:
\[
\begin{array}
{ccc}
{} & {} & Nr(\Z) \\
{} &{}&\Big\downarrow\vcenter{\rlap{$\scriptstyle{\text{int}}$}}\\
\mathbb W(\mathbb Z)
& \stackrel{\tau \, (\cong)}{\longrightarrow}
& \hat \Omega(C)
\end{array}
\quad
\stackrel{\text{if }R \text{ is torsion-free }}{\dashrightarrow}
\quad
\begin{array}
{ccc}
{} & {} & Nr_G(R)\\
{} &{}&\Big\downarrow\vcenter{\rlap{$\scriptstyle{\text{int}}$}}\\
\mathbb W_G(R) & \stackrel{\tau_R\,(\cong) }{\longrightarrow} & \hat \Omega_R(G)
\end{array}
\]

Very recently, developing the above construction,
we have constructed a functor, $\Delta_G$, such that
(i) $\Delta_G$ is equivalent to $\mathbb W_G$, and
(ii) $\Delta_G(\mathbb Z)=\hat \Omega_R(\mathbb Z)$
(\cite{O2}).

On this topic, we first study a general theory of
the Witt-Burnside ring and the necklace ring of a profinite group.
For example, we will show that the $R$-Teichm\"{u}ller map,
$\tau_R:\mathbb W_G(R)\to Nr_G(R)$,
is a ring isomorphism for arbitrary special $\ld$-rings $R$.
Next, we apply the above general theory to the special case
where $G=\hat C$.
Form the fact that $\tau_R$ is a ring homomorphism we obtain that
\begin{equation}\label{isomorphisms}
\mathbb W(R)\cong Nr(R)\cong \Lambda_1(R) \text{ (as rings)  }
\end{equation}
for any special $\ld$-ring $R$.
Here, $\Lambda_1(R)$ is the modified Grothendieck ring of formal power
series with constant term 1.

Grothendieck showed in 1956 that $\Lambda_1(A)$ (originally $\Lambda(A)$)
has a special $\ld$-ring structure for any unital commutative ring $A$.
In view of Eq. \eqref{isomorphisms}, this implies that
if $R$ is a special $\ld$-ring, then
it is possible to make $\mathbb W(R)$ and $Nr(R)$ into
special $\ld$-rings by transporting the special
$\ld$-ring structure of $\Lambda_1(R)$.
Based on this observation, we will show that $\mathbb W$ and $Nr$
can be viewed as functors from the category of special $\ld$-rings
into the category of special $\ld$-rings.
In this paper, we will focus particularly on their Adams
operations rather than their $\ld$-operations since
the former often are much easier to handle than the latter.
Indeed, the Adams operations coincide with the well-known Frobenius operators
(Section \ref{PRELIMINARIES} and Section \ref{Section 3:Necklace ring}).

Next, we study necklace polynomials.
In particular, we investigate the class of rings over which
the necklace polynomials
\begin{equation}\label{necklace 0*}
\mathcal M(x,n):= \dfrac 1n\displaystyle\sum_{d|n} \mu(d)\,x^{\frac
nd}, \quad n\in \mathbb N,
\end{equation}
remain valid.

The second topic relates to the logarithmic function.
The concept of the logarithmic function was first introduced by R. M.
Bryant (\cite {B}) to explain the properties of {\it Lie module
functions} of free Lie algebras.
Let $G$ be a group, $K$ a field, and $V$ a $KG$-module. Letting
$\mathfrak L(V)$ be the free Lie algebra generated by $V,$ $G$
acts on $\mathfrak L(V)$ by Lie algebra automorphisms.
In particular, $\mathfrak L(V)_n,$ the $n$-th homogeneous components
in $\mathfrak L(V),$ are themselves $KG$-modules.
The {\it Lie module function} of
$\mathfrak L(V),$ denoted by $[\mathfrak L(V)],$ is a formal
$q$-series,
$$\sum_{n\ge 1}[\mathfrak L(V)_n] q^n,$$
where the
coefficients $[\mathfrak L(V)_n]$ are the elements of the Green
algebra, $\G_K(G)$ of $G$ over $K$ (Section 5.2).
Bryant proved that for all $G$ and $K,$
there exists a unique logarithmic
function, $\mathcal D$, on $q\G_K(G)[[q]]$, such that
\begin{equation}\label{bryantI}
[\mathfrak L(V)]=\mathcal D([V])
\end{equation}
for every $\mathbb N$-graded $KG$-module $V.$
However, despite of the existence of such a
logarithmic function, we have no explicit idea of how it looks in
general.
Consequently, he introduced the concept of the {\it
Grothendieck Lie module function} of $\mathfrak L(V),$ denoted by
$\overline {\mathfrak L(V)},$ instead of the Lie module function
$[\mathfrak L(V)]$, when $G$ is a finite group.

The {\it Grothendieck Lie module function},
$\overline {\mathfrak L(V)}$, is a formal $q$-series defined by
$$\sum_{n\ge 1}\overline { \mathfrak
L(V)_n}q^n,$$ where the coefficients $\overline {\mathfrak
L(V)_n}$ are the elements of the Grothendieck algebra $\overline
\G_K(G)$ of $G$ over $K$ (Section 5.2).
In this situation, he
not only showed that for every finite group, $G$,
and every field, $K,$
there exists a unique logarithmic function,
$\overline{\mathcal D}$, on $ q\overline {\G}_K(G)[[q]]$, such that
\begin{equation}\label{bryantII}
\overline {\mathfrak L(V)}=\overline {\mathcal D}(\overline V)
\end{equation}
for all $KG$-modules $V,$ but he also provided its explicit form.

In this paper, we generalize Bryant's results to the case of
graded Lie superalgebras with a group action.
Let $\widehat {\Gamma}$ be a free abelian group with finite rank,
and let
$\Gamma$ be a countable (usually infinite) sub-semigroup in
$\widehat {\Gamma}$, such that every element, $\alpha \in \Gamma$,
can be written as a sum of elements of $\Gamma$ in only
a finite number of ways.
Consider ($\G \times \mathbb Z_2$)-graded Lie superalgebras,
$$\mathfrak L=\oplus_{(\alpha,a)\in \G \times \mathbb Z_2}
\mathfrak L_{(\alpha,a)},$$
with $\dim \mathfrak L_{(\alpha,a)} < \infty $ for
all $(\alpha,a) \in (\G\times\mathbb Z_2).$
Suppose $G$ acts on $\mathfrak L$,
preserving the ($\G \times \mathbb Z_2$)-gradation.
Our goal in this work is to find the condition
on $G$ and $K$ for which there exists a unique logarithmic function
on $\G_K(G)[[\G \times \mathbb Z_2]]$
(resp., $\overline \G_K(G)[[\G \times \mathbb Z_2]]$))
satisfying an identity such as \eqref {bryantI} (resp., \eqref {bryantII}).
We also seek to compute the explicit form of this logarithmic function.

This paper is organized as follows.
In Section \ref{PRELIMINARIES}, we introduce the basic
definitions and notation used in this paper.
We also introduce a $q$-deformation of the Grothendieck ring
of formal power series with constant term 1, which is isomorphic to
the $q$-deformation of the universal ring of Witt vectors.
In Section \ref{Section 3:Necklace ring},
we present the main results on the necklace polynomial
and the necklace ring.
In Section \ref{LOGARITHM}, we prove the main theorems of a logarithmic function
(Theorems \ref{logarithm function}, \ref {green ring-logarithm function},
and \ref {char k-logarithm function}) for graded Lie superalgebras.
As corollaries of these theorems, we obtain closed formulas for
the homogeneous components, $|\fL_{(\al,a)}|$ and $|\overline
{\fL_{(\al,a)}}|$ (Corollaries \ref{root space1 }
and \ref{chark-root space1 }).
The final section is devoted to applications.
In Section \ref {PLETHY} we interpret the symmetric
power map, $\tilde s_t$, using plethysm
and then make a few remarks on symmetric functions.
In Section \ref {supersymmetric*}
we present several generating sets of supersymmetric functions.
Section \ref{RECURSIVE FORMULA}
contains recursive formulas for computing $[\fL_{(\al,a)}]$.
Frequently, recursive formulas
are more efficient than closed formulas for this purpose,
so, we will give recursive formulas for
$|\fL_{(\al,a)}|$ and $|\overline {\fL_{(\al,a)}}|$
(Proposition \ref {recursive formula 1}).
Section \ref{REPLICABLE} discusses a new
interpretation of completely replicable functions from a viewpoint of
logarithmic functions and $\ld$-ring structures (Propositions \ref
{characterization OF REP I} and \ref{characterization OF REP II}).

\section{PRELIMINARIES}\label{PRELIMINARIES}

\subsection{$\ld$-rings and Grothendieck
ring of formal power series}\label{preliminary:1}
A {\it $\ld$-ring} $R$ is a unital commutative ring with operations
$\ld^n : R \to R, \,
(n=0,1,2,\cdots)$ such that
\begin{align*}
& (1)\quad \ld^0(x)=1,\\
& (2)\quad \ld^1(x)=x,\\
& (3)\quad \ld^n(x+y)=\sum_{r=0}^{n}\ld^r(x)\ld^{n-r}(y).
\end{align*}
If $t$ is an indeterminate, we let
$$\ld_t(x):=\sum_{n=0}^{\infty}\ld^n(x)t^n, \quad x \in R.$$
By the third condition, it is
straightforward that $\ld_t(x+y)=\ld_t(x)\ld_t(y).$
An element $x\in R$ is said to be {\it $n$-dimensional}
if $\ld_t(x)$ is a polynomial of degree $n$ in $t$.
For further information refer to \cite{AD,H,DK}.

Grothendieck, in 1956, introduced a functor, $\Lambda$, from the
category of unital commutative rings to the category of
special $\ld$-rings.
$\Lambda(A)$ is the ring whose underlying set is
$$1+A[[t]]^+:=\{ 1+\sum_{n=1}^{\infty} x_nt^n  \,:\, \,x_n \in A\}.$$
Its ring structure is determined in the following way:
Let $x_1,x_2,\cdots ;\,
y_1,y_2,\cdots$ be indeterminates.
We define $s_i $ (resp.
$\sigma_j$) to be the elementary symmetric functions in variables
$x_1,x_2,\cdots$ (resp. $y_1,y_2,\cdots$), that is,
\begin{align*}
&(1+s_1 t +s_2t^2+\cdots)=\prod_{i=1}^{\infty}(1+x_i t),\\
&(1+\sigma_1 t +\sigma_2t^2+\cdots)=\prod_{i=1}^{\infty}(1+y_i t).
\end{align*}
Set $P_n(s_1,\cdots,s_n;\sigma_1,\cdots,\sigma_n)$ to be the
coefficient of $t^n$ in $$\prod _{i,j \ge 1}(1+x_iy_jt),$$ and
$P_{n,m}(s_1,\cdots,s_{mn})$ the coefficient of $t^n$ in
$$\prod_{i_1<i_2<\cdots <i_m}(1+x_{i_1}\cdots x_{i_m}t).$$
With this notation, let us define the $\ld$-ring structure on $\Lambda(A)$ by

(a) $\oplus$ : Addition is just the multiplication of power series.\\
(b) ${\star} $ : Multiplication is given by
$$
(1+\sum a_nt^n)\star (1+\sum b_nt^n)=1+\sum
P_n(a_1,\cdots,a_n;b_1,\cdots,b_n)t^n.$$
(c) $ \Lambda^m(1+\sum a_nt^n)=1+\sum
P_{n,m}(a_1,\cdots,a_{nm})t^n.$

The operations of $\Lambda(A)$
can be understood better in terms of symmetric functions.
Let $X=\{x_i:i\ge 1 \}$ and $Y=\{y_j: j\ge 1 \}$ be
the infinite sets of commuting indeterminates $x_i$'s and $y_i$'s respectively.
We call these sets {\it alphabets}.
Now, let us introduce the following operations on alphabets:
\begin{align*}
&X+Y=\{x_i,y_i: i\ge 1\},\\
&X\cdot Y=\{x_iy_j: i,j\ge 1\},\\
&\Lambda^m(X)=\{x_{i_1}\cdots x_{i_m}: i_1<i_2<\cdots <i_m\}.
\end{align*}
Exploiting the notation
$$E(X,t)=\prod_{i=1}^{\infty}(1+x_i t),$$
then the above conditions (a) through (c)
can be regarded as expressing the identities
\begin{align*}
&E(X,t)\oplus E(Y,t) =E(X+Y,t),\\
&E(X,t)\star E(Y,t)=E(X\cdot Y,t),\\
&\Lambda^m(E(X,t))=E(\Lambda^m(X),t).
\end{align*}
Similarly, we can introduce another ring structures on $1+A[[t]]^+$
via the bijective maps
\begin{align*}
&\theta_0:\Lambda(A)\to 1+A[[t]]^+, \quad f(t)\mapsto f(-t),\\
&\theta_1:\Lambda(A)\to 1+A[[t]]^+, \quad f(t)\mapsto
\frac{1}{f(-t)}.
\end{align*}
Denote by $\Lambda_i(A)\,(i=0,1)$ the ring whose ring structure is
induced from $\theta_i,\,(i=0,1).$
In this paper, we mainly deal with the ring $\Lambda_1(A)$ rather than
$\Lambda(A)$ and $\Lambda_0(A)$.
Let $\oplus$ (resp., $\star_1$) denote the addition (resp., the multiplication)
of $\Lambda_1(A)$. Define $\bar s_i $ (resp., $\bar \sigma_j$) to
be the symmetric functions in variables $x_1,x_2,\cdots$ (resp.,
$y_1,y_2,\cdots$) determined by the equations
\begin{align*}
&(1+\bar s_1 t + \bar s_2t^2+\cdots)=\prod_{i=1}^{\infty}\frac {1}{1-x_i t},\\
&(1+\bar \sigma_1 t +\bar
\sigma_2t^2+\cdots)=\prod_{i=1}^{\infty}\frac {1}{1-y_i t}.
\end{align*}
Set $\bar P_n(\bar s_1,\cdots,\bar s_n;\bar \sigma_1,\cdots,\bar
\sigma_n)$ to be the coefficient of $t^n$ in $$\prod_{i,j}\frac
{1}{1-x_iy_jt},$$ and $\bar P_{n,m}(\bar s_1,\cdots,\bar s_{mn})$
the coefficient of $t^n$ in
$$\prod_{i_1<i_2<\cdots <i_m}\frac{1}{1-x_{i_1}\cdots x_{i_m}t}.$$
From the fundamental theorem on symmetric functions (see \cite {Mac})
it follows that
\begin{align*}
&\bar P_n(X_1,\cdots,X_n;Y_1,\cdots,Y_n)\in \mathbb Z[X_1,\cdots,X_n;Y_1,\cdots,Y_n]\\
&\bar P_{n,m}(X_1,\cdots,X_{mn})\in \mathbb Z[X_1,\cdots,X_{mn}].
\end{align*}
Indeed, the $\ld$-ring structure on $\Lambda_1(A)$ is given by

(a') $\oplus$ : Addition is just multiplication of power series.\\
(b') ${\star}_1 $ : Multiplication is given by
$$
(1+\sum a_nt^n)\star_1 (1+\sum b_nt^n)=1+\sum \bar
P_n(a_1,\cdots,a_n;b_1,\cdots,b_n)t^n.$$
(c') $ \bar \Lambda^m(1+\sum a_nt^n)=1+\sum \bar
P_{n,m}(a_1,\cdots,a_{nm})t^n.$

As in the previous paragragh, if we let
$$H(X,t)=\prod_{i=1}^{\infty}\frac{1}{1-x_i t},$$
then the above conditions (a') through (c')
can be regarded as expressing the identities
\begin{align*}
&H(X,t)\oplus H(Y,t)=H(X+Y,t),\\
&H(X,t)\star_1 H(Y,t)=H(X\cdot Y,t),\\
&\bar \Lambda^m(H(X,t))=H(\Lambda^m(X),t).
\end{align*}

\begin{df}
A $\ld$-ring $R$ is said to be {\it special} if $\ld_{t}: R \to
\Lambda(R)$ is a $\ld$-homomorphism, that is, a ring homomorphism
commuting with the $\ld$-operations.
\end{df}

\subsection{Adams operations and binomial rings}\label{preliminary:2}
Let $R$ be a $\ld$-ring.
We define the {\it $n$-th Adams operation, $\Psi^n :R\to R$}, by
\begin{equation}\label {adams}
\frac{d}{d{\it t}}\log
\ld_t(x)=:\sum_{n=0}^{\infty}(-1)^n\Psi^{n+1}(x)t^n
\end{equation}
for all $x \in R.$
If $R$ is a special $\ld$-ring, it is well-known that
$\Psi^n$ are $\ld$-ring homomorphisms and $\Psi^n \circ \Psi^m=\Psi^{mn}$
for all $m,n\ge 1.$
On the other hand, if we define the {\it $n$-th symmetric power operations} by
\begin{equation}\label{symmetric map}
\mathfrak S_t(x)=\sum_{n=0}^{\infty}S^n(x)t^n:=\dfrac{1}{\ld_{-t}(x)},
\end{equation}
then Eq. \eqref {adams} can be rewritten as
\begin{equation}\label {symm power operation}
\frac{d}{d{\it t}}\log
\mathfrak S_t(x)=\sum_{n=0}^{\infty}\Psi^{n+1}(x)t^n.
\end{equation}
With this notation, it is easy to show that $R$ is special if and only if
$\mathfrak S_t: R \to \Lambda_1(R)$ is a $\ld$-homomorphism.

Frequently, Adams operations completely determine the structure of
the associated $\ld$-ring.
A ring $R$ is said to be a {\it $\Psi$-ring} if it is a unital
commutative ring with a set of operations $\Psi^n: R\to
R$ for all $n \ge 1,$ satisfying
\begin{equation*}
\begin{aligned}
&\Psi^1(a)=a,\\
&\Psi^n(a+b)=\Psi^n(a)+\Psi^n(b)
\end{aligned}
\end{equation*}
for all $a,b \in R.$ A {\it special $\Psi$-ring} is defined to be
a $\Psi$-ring satisfying
\begin{equation*}
\begin{aligned}
&\Psi^n(1)=1,\\
&\Psi^n(ab)=\Psi^n(a)\Psi^n(b),\\
&\Psi^n(\Psi^m(a))=\Psi^{nm}(a)
\end{aligned}
\end{equation*}
for all $m,n\ge 1$ and all $a,b \in R.$
Note that if $R$ is a unital commutative ring, then it becomes
a special $\Psi$-ring by setting $\Psi^n=id$ for all $n\ge 1$.
\begin{thm}{\rm (\cite{W})}\label {psi=ld}
Let $R$ be a special $\Psi$-ring which has no $\mathbb Z$-torsion
and such that $\Psi^p(a)=a^p$ mod $pR$ if $p$ is a prime.
Then, there is a unique special $\ld$-ring structure on $R$ such that
the $\Psi^n$ are the associated Adams operations.
\end{thm}
As an easy application of Theorem \ref {psi=ld}, we have
\begin{cor}\label{psiring}
Let $R$ be a commutative ring with unity which has no $\mathbb Z$-torsion
and such that $a^p=a$ mod $pR$ if $p$ is a prime.
Then, $R$ has a unique special $\ld$-ring structure
with $\Psi^n=id$ for all $n\ge 1$.
\end{cor}
\noindent{\bf Proof.}
For all $n\in \mathbb N$ let $\Psi^n$ be the identity map on $R$.
Then, it is clear that $R$ is a special $\Psi$-ring with regard to the operators
$\{\Psi^n:n\in \mathbb N\}$. Now, Theorem \ref {psi=ld} implies our assertion.
\qed

\begin{exm}
Let $R$ be a $\mathbb Q$-algebra, $\Z$, or $\mathbb Z_{(r)}$
the ring of integers localized at $r$.
Then one can easily verify that it satisfies the condition of Corollary
\ref{psiring}.
\end{exm}

\begin{df}
A special $\ld$-ring in which $\Psi^n={\rm id}$ for all $n\ge 1$
will be called a {\it binomial ring}.
\end{df}

\subsection{$q$-deformation of Grothendieck ring of formal power series}
It is quite interesting to note that the universal ring of Witt vectors
has a $q$-analogue for every integer $q$.
Recall that the universal ring $\mathbb W(A)$ of Witt vectors over $A$
is isomorphic to $\Lambda_1(A)$.
Thus, it would be natural to think over the q-analogue of $\Lambda_1(A)$
which corresponds to the $q$-deformation of $\mathbb W(A)$
(refer to \cite{O2}).
In this section, we deal with the question very briefly.

We start with remarking that notations associated with
formal group laws used in this paper can be found in \cite{L}.
In the theory of formal group laws, it is well-known that
$\Lambda_1(A)$ is isomorphic (as abelian groups) to the group of curves,
$$\mathcal C(F_1,A)=\{\sum_{n\ge 1}a_nt^n: a_n\in A \},$$
in the multiplicative formal group law
$F_1(X,Y)=X+Y-XY$
via the isomorphism
$$\beta: \gamma(t)\mapsto \frac{1}{1-\gamma(t)}.$$

More generally, we can verify that the group of curves $\mathcal C(F_q,A)$
in the formal group law
$F_q(X,Y)=X+Y-qXY,\,(q\in \mathbb Z\setminus \{0\}),$
is isomorphic (as abelian groups) to $\Lambda_1(A)$ via the isomorphism
$$\beta^q: \gamma(t)\mapsto \frac{1}{1-q\gamma(t)}$$
if $A$ is a ring in which $q$ is invertible.
On the other hand, $\mathcal C(F_q,A)$ is isomorphic to
$\mathbb W^{F_q}(A)$, the universal ring of Witt vectors over $A$
associated with $F_q$,
via the Artin-Hasse type exponential map
$$H^{F_q}:\mathbb W^{F_q}\to \mathcal C(F_q,A),
\quad \alpha \mapsto {\sum_{n\ge 1}}^{F_q}\alpha_n t^n.$$
Here, it should be remarked that
if we endow $\mathcal C(F_q,A)$ with
the ring structure via $H^{F_q}$
then $\beta^q$ is no longer a ring homomorphism.
In order to make it into a ring
homomorphism, we need to modify the multiplication of
$\Lambda_1(A)$ so as to satisfy
\begin{equation}\label{q-deform of Grothendieck}
\prod_{i=1}^{\infty}\left(\frac
{1}{1-x_i t}\right)^q \star_q\prod_{j=1}^{\infty}\left(\frac {1}{1-y_j t}\right)^q
=\prod_{i,j}\left(\frac {1}{1-x_iy_j t}\right)^q.
\end{equation}
Here, the notation $\star_q$ represents the modified multiplication.
This can be verified by observing that
the map, $E^{F_q}: {\rm gh}(A)\to \mathcal C(F_q,A)$, is given by
$$(\beta^q)^{-1}\circ \exp\circ \iota_q$$
since $$\log_q(X)=\sum_{n\ge 1}\frac {q^{n-1}}{n}X^n=
\frac 1q\log \left(\frac {1}{1-qX}\right).$$
Here, the map $\iota_q:{\rm gh}(A)\to A[[t]]$ is defined by
$$(a_n)_n \mapsto \sum_{n\ge 1}\frac 1n a_nt^n.$$
Combining the relation
$$(\beta^q\circ E^{F_q})^{-1}\left(\prod_{i=1}^{\infty}\left(\frac
{1}{1-x_i t}\right)^q\right)=(p_n(X))_n$$
with the well-known identity
$$p_n(X)\cdot p_n(Y)=p_n(X\cdot Y)$$
yields the identity \eqref{q-deform of Grothendieck}.
The notation $p_n(X)$ represents the $n$-th power sum of the alphabet $X$.
From this it follows that
$$\prod_{i=1}^{\infty}\left(\frac{1}{1-x_i t}\right)^q
\star_q \prod_{i=1}^{\infty}\left(\frac {1}{1-y_i t}\right)^q
=(\beta^q\circ E^{G^q_m})((p_n(X\cdot Y))_{n\ge 1}).$$ Denote by
$\Lambda_q(A)$ the ring induced from Eq. \eqref{q-deform of
Grothendieck}. In conclusion, if $A$ is a commutative ring in which $q$ is invertible,
then
$$\mathbb W^{F_q}(A)\cong \mathcal C(F_q,A)\cong \Lambda_q(A)\quad (\text{as rings}).$$
\subsection{Green algebras and Grothendieck algebras}
We end this section by recalling some terminologies which will be
exploited in the part of logarithmic functions.

Let $G$ be a group and $K$ be a field.
We consider a set $\{I_\lambda \,:\,\lambda \in \Lambda \}$
consisting of representatives from each isomorphism class
of finite dimensional indecomposable $KG$-modules.
If $V$ is any finite dimensional $KG$-module,
then we write $[V]$ for the element of $\sum\alpha_\lambda I_\lambda,$
where $\alpha_\lambda$ is the number of summands isomorphic to
$I_\lambda$ in an unrefinable direct sum decomposition of $V.$
The free abelian group generated by $\{ I_\lambda \,:\,\lambda \in
\Lambda \},$ denoted by $R_K(G),$
is called the {\it Green ring} (or {\it representation ring})
{\it of $G$ over $K,$} where multiplication on $R_K(G)$
is defined by tensor products.
Note that for all finite dimensional $KG$-modules
$U$ and $V,$ $[U]+[V]=[U\oplus V]$ and $[U][V]=[U\otimes V].$
The scalar extension $\C \otimes_{\mathbb Z} R_K(G),$
which becomes a commutative $\mathbb C$-algebra,
is called the {\it Green algebra of $G$ over $K$}
and denoted by $\Gamma_K(G).$ In particular, when
$K=\C,$ we use the notations $R(G),\,\Gamma(G)$ instead of
$R_K(G),\, \Gamma_K(G).$

From now on, let $G$ be a finite group and $K$ be an arbitrary
field. Let $I$ be the subspace of $\Gamma_K(G)$ spanned by all the
elements of the form $[V]-[U]-[W],$ where $U,V,W$ are finite
dimensional $KG$-modules occurring in a short exact sequence $0\to
U \to V \to W\to 0.$ Note that $I$ is an ideal. We consider
$\overline \Gamma_K(G)=\Gamma_K(G)/I$, the quotient of
$\Gamma_K(G)$ by $I$, which is called the {\it Grothendieck
algebra}. Similarly we define the {\it Grothendieck group}
$\overline R_K(G).$ Then, we have $\overline \G_K(G)= \C \otimes
\overline R_K(G).$ If $V$ is any finite dimensional $KG$-module,
then we write $\overline V$ for the element of $\overline \G_K(G)$
and $\overline R_K(G).$
\section{NECKLACE RINGS}\label{Section 3:Necklace ring}

\subsection{The Witt-Burnside ring
and the necklace ring of a profinite group}
\label{sec3.1}
In this section, we start with recalling a functor denoted by $Nr_G$.
This functor was first introduced as a generalization
of the necklace ring functor of Metropolis and Rota (see \cite{MR,O}).
In the next section, we focus on the case where $G$ is $\hat C$,
the profinite completion of the multiplicative infinite cyclic group $C$.

Let $G$ be a profinite group.
For two subgroups $U,\,V$ of $G$, we say that
$U$ is subconjugate to $V$ if $U$ is a subgroup
of some conjugates of $V$.
This is a partial order on the set of
the conjugacy classes of open subgroups
of $G$, and will be denoted by $[V]\preceq [U]$.
Consider (and fix) an enumeration of this poset
satisfying the condition
$$[V]\preceq [U] \Rightarrow
[V] \text{ precedes }[U].$$
By abuse of notation denote this poset by $\mathcal O(G)$.

Let $R$ be a special $\ld$-ring.
With the above notation, let us define
the {\it necklace ring of $G$ over $R$}, denoted by $Nr_G(R)$, by
the ring whose underlying set is
$$\prod_{[U]\in \mathcal O(G)}R.$$
Its addition is defined componentwise.
On the contrary, its multiplication is somewhat complicated.
For $x=(x_{[V]})_{[V]}$ and $y=(y_{[W]})_{[W]}$ define
\begin{equation*}
(x\cdot y)_{[U]}:=\sum_{[V],[W]}\sum_{VgW\subseteq G \atop [Z(g,V,W)]=[U]}
\Psi^{(V:Z(g,V,W))}(x_{[V]})\Psi^{(W:Z(g,V,W))}(y_{[W]}),
\end{equation*}
Here the notation $Z(g,V,W)$ represents $V\cap gWg^{-1}$
and the notation $(V:U)$ means the index of $U$ in $gVg^{-1}$.

For an open subgroup $U$ of $G$
let us introduce the induction map,
\begin{equation}\label{defn of inductions}
\text{Ind}_U^G(R):Nr_U(R)\to Nr_G(R),\quad
(x_{[V]})_{[V]\in \mathcal O(U)}
\mapsto  (y_{[W]})_{[W] \in \mathcal O(G)},
\end{equation}
where
$$y_{[W]}=\sum_{[V]=[W] \text{ in } \mathcal O(G)}x_{[V]}.$$
Also, we introduce the exponential map
\begin{equation*}
\tau_R^G: R \to Nr_G(R),\quad r\mapsto
\left(\sum_{[V]\in \mathcal O(G)}M_G(r,V)\right)_{[V]}.
\end{equation*}
Here the coefficients $M_G(r,V)$ are determined in the following manner.
First, we write $r$ as a sum of one-dimensional elements, say,
$r_1+r_2+\cdots+r_m.$
Now, consider $\mathcal C(G,{\mathfrak r})$
the set of continuous maps from $G$ to the topological space
$${\mathfrak r}:=\{r_1,r_2,\cdots,r_m\}$$
with regard to the discrete topology with trivial $G$-action.
It is well-known that $\mathcal C(G,{\mathfrak r})$
becomes a $G$-space
with regard to the compact-open topology
via the following standard $G$-action
$$(g\cdot f)(x)=f(g^{-1}\cdot x)$$
(see \cite {DS2}).
Write $\mathcal C(G,{\mathfrak r})$ as the disjoint union of $G$-orbits,
say,
$$\underset{h}{\dot\bigcup}
\,G\cdot h,$$
where $h$ runs through a system of representatives of this decomposition.
After writing $G/G_h
=\underset{1\le i \le (G:G_h)}{\dot\bigcup}w_iG_h$,
where $G_h$ represents the isotropy subgroup of $h$,
we let
$$[h]:=\prod_{i=1}^{ (G:G_h)}h(w_i).$$
Clearly, this is well-defined since $h$ is $G_h$-invariant.
With this notation, we define $M_G(r,V)$ by
$$\sum_h [h],$$
where $h$ is taken over the representatives such that
$Gh$ is isomorphic to $G/V$.

Recently we have constructed a map, called
{\it $R$-Teichm\"{u}ller map},
\begin{align*}
\tau_R: \mathbb W_G(R)\to Nr_G(R),\quad
\alpha\mapsto \sum_{[U]\in \mathcal O(G)}\,\text{Ind}_U^G
(\tau_R^U(\alpha([U]))),
\end{align*}
and then have shown that this map is bijective, and a ring homomorphism
if $R$ is torsion-free (see \cite{O}).

\begin{rem}
Note that we have identified $Nr_G(R)$ with $\hat \Omega_G(R)$
via the interpretation map {\rm (}for short, int{\rm )} which is given by
$$\text{\rm interpretation}:\hat \Omega_R (G)\to Nr_G(R),
\quad \sum_{[U]}a_{[U]}[G/U]\mapsto (a_{[U]})_{[U]}.$$
\end{rem}

The following lemma is immediate by
the definition of the isomorphism $\tau_R$.
\begin{lem}\label {def of teichmuller}
Let $[W]\in \mathcal O(G)$.
Then, the $[W]$-th component of $\tau_R(\alpha), \alpha \in \mathbb W_G(R),$
is given by
\begin{equation*}
\sum_{[U]\in \mathcal O(G)}
\sum_{[V]\in \mathcal O(U)}M_U(\alpha([U]),V),
\end{equation*}
where $[V]$ ranges over the elements $[Z] \in \mathcal O(U)$
such that $[Z]=[W]$ in $\mathcal O(G)$.
\end{lem}

\begin{thm}\label{ringhomo of teichmuller}
$\tau_R$ is a ring isomorphism.
\end{thm}

\noindent{\bf Proof.}
It was shown in \cite {O} that $\tau_R$ is bijective
for every special $\ld$-ring.
So, for our purpose we have only to show that it is a ring homomorphism.
Given a unital commutative ring $R$,
Dress and Siebeneicher \cite {DS2} showed that
for $\alpha, \beta \in \mathbb W_G(R)$
$$\alpha+\beta
=(s_U(\alpha_{[V]},\beta_{[V]}\,|\,[G]\preceq[V]\preceq [U]))_{[U]\in \mathcal O(G)}$$
and
$$\alpha \cdot\beta
=(p_U(\alpha_{[V]},\beta_{[V]}\,|\,[G]\preceq[V]\preceq [U]))_{[U]\in \mathcal O(G)}$$
for some integral polynomials $s_U$ and $p_U$
for every
$[U]\in \mathcal O(G)$.
Observe that for every $[Z]\in \mathcal O(G)$,
the $[Z]$-th component of
$\tau_R({\bf x}), \,{\bf x}=(x_{[U]})_{[U]\in \mathcal O(G)}$,
is an integral polynomial in $\ld^k(x_{[U]})$'s for
$$1\le k \le (U:Z),\,\, [G]\preceq[U]\preceq [Z].$$
This follows from Lemma \ref {def of teichmuller}.
More precisely, this is because a necklace polynomial
$M_U(\alpha([U]),Z)$ becomes a symmetric polynomial
in one-dimensional elements after writing $\alpha([U])$
as a sum of one-dimensional elements.
Applying the elementary theory of symmetric functions, one can show that
$$M_U(\alpha([U]),Z)$$ can be expressed as an integral polynomial in
$\ld^k(\alpha([U]))$'s, $1\le k \le (U:Z).$
Combining these two observations,
we can conclude that the identities
\begin{equation}\label{universal identity}
\begin{aligned}
&\tau_R(\alpha+\beta)=\tau_R(\alpha)+\tau_R(\beta),\\
&\tau_R(\alpha \cdot \beta)=\tau_R(\alpha)\cdot \tau_R(\beta)
\end{aligned}
\end{equation}
hold universally.
In other words, two equations in \eqref{universal identity}
are always true for arbitrary special  $\ld$-rings.
This completes the proof.
\qed

As a byproduct of Theorem \ref{ringhomo of teichmuller}
we obtain an trivial, but very significant corollary.
Suppose that $R$ has two kinds of special $\ld$-ring structures.
Then, each structure produces different necklace rings,
say, $Nr_{G,1}$ and $Nr_{G,2}$.
\begin{cor}\label{canonical isom}
$Nr_{G,1}(R)$ is canonically isomorphic to $Nr_{G,2}(R)$.
\end{cor}
\noindent{\bf Proof.}
Consider the map
\begin{equation*}
Nr_G({\rm id}_{12}):= \tau_{R,2} \circ{\tau_{R,1}}^{-1}:Nr_{G,1}(R)\to Nr_{G,2}(R).
\end{equation*}
Clearly $Nr_G({\rm id}_{12})$ is a ring-isomorphism
by Theorem \ref{ringhomo of teichmuller}. \qed

From now on, we investigate inductions and restrictions on
$$\hat \Omega_R(G),\quad \mathbb W_G(R),\,\,\text{ and } R^{\mathcal O(G)}$$
which have played a crucial role in in the theory of
necklace rings and Witt-Burnside rings.
In the next section, we show that
if $G=\hat C$ then restriction maps on
these rings coincide with Adams operations associated with
the special $\ld$-structure induced from that of $\Lambda_1(R)$.

First, let us recall inductions and restrictions on
$\hat \Omega_R(G)$, equivalently on $Nr_G(R)$ (see \cite {O,O1}).
For an open subgroup $U$ of $G$,
the induction $\text{Ind}_U^G$,
as defined in Eq. \eqref{defn of inductions},
is an additive homomorphism
from $\hat \Omega_R (U)$ to $\hat \Omega_R(G)$
given by
$$\text{Ind}_U^G\left(\sum_{[W]\in \mathcal O(U)}b_{[W]}[U/W]\right)
=\sum_{[V]\in \mathcal O(G)}\left(\sum_{[W]\in \mathcal O(U)
\atop [W]=[V] \text { in } \mathcal O(G)}b_{[W]}\right) [G/V].$$
While, the restriction
$$\text{\rm Res}_U^G :\hat \Omega_R (G)\to \hat \Omega_R (U)$$
is defined by the rule
\begin{equation*}
\sum_{[V]}b_{[V]}[G/V]\mapsto \sum_{[V]}
\sum_g \Psi^{(V:Z(g,U,V))}
=\sum_{[W]\in \mathcal O(U)}
\left(\sum_{[V]\in \mathcal O(G)\atop [Z(g,U,V)]=[W] \text{ in }\mathcal O(U)}
b_{[V]}\right) [ U/W].
\end{equation*}
Here, $g$ ranges over a set
of representatives of $U$-orbits of $G/V.$
One of many significant properties of restrictions is that
they are indeed ring homomorphisms.
Also, it is worthwhile noting that
for open subgroups $U\leqslant V \leqslant G,$
\begin{equation}\label{transitivity of ins and res}
\begin{aligned}
&\text{\rm Ind}_V^G \circ \text{\rm Ind} _U^V=\text{\rm Ind} _U^G,\\
&\text{\rm Res}_U^V \circ \text{\rm Res}_V^G =\text{\rm Res}_U^G.
\end{aligned}
\end{equation}

In \cite[Lemma 3.13]{O} it has been shown that
the diagram
\begin{equation}\label{commuting diagram}
\begin{picture}(360,70)
\put(100,60){$\mathbb W_G(R)$}
\put(135,63){\vector(1,0){30}}
\put(170,60){$\hat \Omega_G(R)$}
\put(170,6){$R^{\underline {\mathcal O(G)}}$}
\put(175,57){\vector(0,-1){38}}
\put(120,55){\vector(1,-1){45}}
\put(145,67){$\tau_R$}
\put(178,35){$\tilde \varphi$}
\put(125,35){$\Phi$}
\put(155,40){$\curvearrowright$}
\end{picture}
\end{equation}
is commutative.
Here,
$$\Phi(\alpha)=\left(\sum_{[G]\preceq [V]\preceq [U]}
\varphi_U(G/V)\cdot\alpha([V])^{(V:U)}\right)_{[U]}$$
and
\begin{equation*}
\tilde \varphi_U \left(\sum_{[V]\in \mathcal O(G)}b_{[V]}[G/V]\right)=
\left(\sum_{[G]\preceq [V]\preceq [U]}
\varphi_U(G/V)\Psi^{(V:U)}(b_{[V]})\right)_{[U]}\,\,.
\end{equation*}
The notation $\varphi_U(X)$ means the cardinality
of the set $X^U$ of $U$-invariant
elements of $X$ and
let $G/U$ denote the $G$-space of left cosets of $U$ in $G$.
Now, by the transport of inductions and restrictions on $\hat \Omega_G(R)$
via the map $\tau_R$,
we will define the operators ${\bar v}_U,\, {\bar f}_U$ on $\mathbb W_G(R)$.
Indeed, it was already shown in \cite {DS2} that
for an open subgroup $U$ of $G,$
one has well-defined natural transformations
$v_U: \mathbb W_U(-) \to \mathbb  W_G(-)$ and $f_U: \mathbb W_G(-) \to  \mathbb W_U(-)$
satisfying all relations which are known to hold generally for
the restriction
$\text{res}_U^G:\hat \Omega(G)\to \hat \Omega(U)$
and the induction $\text{ind}_U^G:\hat \Omega(U)\to \hat \Omega(G)$.

\begin{thm}\label{natural transformation}
Let $U$ be an open subgroup of $G$.
Regard $\mathbb W_G(-)$ and $\mathbb W_U(-)$ as the functors
from the category of special $\ld$-rings
to the category of commutative rings with identity.
Then, as a natural transformation from $\mathbb W_G(-)$ to $\mathbb W_U(-)$,
${\bar f}_U$ coincides with $f_U$.
Similarly, ${\bar v}_U$ coincides with $v_U$.
\end{thm}
Actually, the proof of Theorem \ref{natural transformation}
can be done essentially in the same way as in \cite{DS2}.
It is based on the following lemmas.
\begin{lem}{\rm (cf. \cite[Lemma (2.12.12)]{DS2})}\label{Lem0}
For any two open subgroups $U,V \leqslant G$ and $\alpha \in R$
one has
\begin{equation*}
\tilde \varphi_U (\,\text{\rm Ind}_V^G(\tau_R^{V}(\alpha))\,)
=\tilde \varphi_U (G/V) \,\alpha^{(V:U)}.
\end{equation*}
\end{lem}
\noindent{\bf Proof.}
Note that
\begin{align*}
&\tilde \varphi_U(\,\text{\rm Ind}_V^G(\tau_R^{V}(\alpha))\,)\\
&=\sum_{gV\in (G/V)^U}\tilde \varphi_U \circ
\text{\rm Res}_U^V(g)(\tau_R^{V}(\alpha))
\text{ (by \cite[Proposition 3.10 (c)]{O})}\\
&=\sum_{gV\in (G/V)^U}\tilde \varphi_U
(\tau_R^{U}(\alpha^{(V\,:\,U)}))
\text{ (by \cite [Lemma 3.11]{O})}\\
&=\tilde \varphi_U(G/V)\,
\alpha^{(V\,:\,U)}.
\end{align*}
\qed

\begin{lem}{\rm (cf. \cite[Lemma (3.2.2)]{DS2})}\label{Lem1}
With the notation in {\rm \cite[Lemma (3.2.2)]{DS2}},
we obtain that for any $\alpha,\beta \in R$
\begin{equation}\label{expasion of addition}
\tau_R^{G}(\alpha+\beta)
=\sum_{G\cdot A \in G\setminus \mathfrak U(G)}
\text{\rm Ind}_{U_A}^{G}(\tau_R^{U_A}(\alpha^{i_A}\cdot \beta^{i_{G-A}})).
\end{equation}
\end{lem}
\noindent{\bf Proof.}
First, we assume that $R$ is torsion-free.
For all open subgroups $U\leqslant G$,
if we take $\tilde \varphi_U$ on the right side of
Eq. \eqref{expasion of addition}, one has
\begin{align*}
&\tilde \varphi_U(\sum_{G\cdot A \in G\setminus \mathfrak U(G)}
\text{\rm Ind}_{U_A}^{G}(\tau_R^{U_A}(\alpha^{i_A}\cdot \beta^{i_{G-A}})))\\
&=\sum_{G\cdot A \in G\setminus \mathfrak U(G)\atop U\lesssim U_A}
\tilde \varphi_U(G/U_A)(\alpha^{i_A}\cdot \beta^{i_{G-A}})^{(U_A:U)}
\text{ (by Lemma \ref{Lem0})}\\
&=\sum_{A \in \mathfrak U(G),\, U \lesssim U_A}
\alpha^{\sharp (A/U)}\cdot \beta^{\sharp (G-A)/U}\\
&=(\alpha+\beta)^{\sharp (G/U)}\\
&=\tilde \varphi_U(\tau_R^G(\alpha+\beta)).
\end{align*}
Since $\tilde \varphi_U$ is injective, we have the desired result.
However, in case where $R$ is not torsion-free,
$\tilde \varphi_U$ is no longer injective.
In this case, we note that the $[U]$-th components appearing in
both sides of Eq. \eqref{expasion of addition} are integral polynomials in
$\ld^k(\alpha)$'s and $\ld^l(\beta)$'s for
$1 \le k,l \le  [G:U]$.
This implies that Identity \eqref{expasion of addition}
holds regardless of torsion.
\qed

\begin{lem}{\rm (cf. \cite[Lemma (3.2.5)]{DS2})} \label{Lem2}
For some $k\in \mathbb N$ let $V_1,\cdots,V_k \leqslant G$
be a sequence of open subgroups of $G$.
Then, for every open subgroups $U \leqslant G$ and every sequence
$\varepsilon_1,\cdots,\varepsilon_k \in \{\pm 1\}$
there exists a unique polynomial
$\xi_U
=\xi^G_{(U;V_1,\cdots,V_k;\varepsilon_1,\cdots,\varepsilon_k)}
=\xi_U(x_1,\cdots,x_k)\in \mathbb Z[x_1,\cdots,x_k]$
such that for all $\alpha_1,\cdots,\alpha_k \in R$
one has
$$\tau^{-1}_R(\sum_{i=1}^k
\varepsilon_i\cdot \text{\rm Ind}_{V_i}^G(\tau_R^{V_i}(\alpha_i))(U)
=\xi_U(\alpha_1,\cdots,\alpha_k).$$
\end{lem}

\noindent{\bf Proof.}
The proof can be done in the exactly same way of that of
\cite[Lemma (3.2.5)]{DS2}.
\qed

Note that the polynomial $\xi_U(x_1,\cdots,x_k)$ must coincide with
that of \cite[Lemma (3.2.5)]{DS2} since $R$ contains $\mathbb Z$
and
$$\tau_R=\tau,\quad \text{Ind}_U^G
=\text{ind}_U^G \quad \text{ if }R=\mathbb Z.$$

\noindent{\bf Proof of Theorem \ref{natural transformation}.}
For $\alpha \in \mathbb W_G(R)$ one has
\begin{align*}
&\text{Res}_U^G(\tau_R(\alpha))\\
&=\sum_{[V]}\text{Res}_U^G\cdot \text{Ind}_V^G(\tau^V_R(\alpha([V])))\\
&=\sum_{[V]}\sum_{UgV\subseteq G}
\text{Ind}_{U\,\cap \, gVg^{-1}}^U \cdot
\text{Res}_{U \, \cap \, gVG^{-1}}^V(g)(\tau^V_R(\alpha([V])))\\
&=\sum_{[V]}\sum_{UgV\subseteq G}
\text{Ind}_{U\,\cap \, gVg^{-1}}^U
(\tau_R^{U \, \cap \,gVg^{-1}}(\alpha([V])^{(V:U \, \cap \,gVg^{-1})}))\\
\end{align*}
Hence, with the same notation as in \cite [Eq. (3.3.9)]{DS2},
it follows from Lemma \ref{Lem2} that
for any open subgroup $W$ of $U$
\begin{equation}
\begin{aligned}
&(\tau_R^{-1}\circ \text{Res}_U^G\circ \tau_R)(\alpha)([W])\\
&=\xi^U_{(W;W_1,\cdots,W_k;1;\cdots;1)}
(\alpha([V_1])^{(V_1;W_1)},\cdots,\alpha([V_k])^{(V_k,W_k)}).
\end{aligned}
\end{equation}
Similarly, one can show
$(\tau_R^{-1}\circ \text{Ind}_U^G\circ \tau_R)(\alpha)([W])$
is a polynomial with integral coefficients in those $\alpha([V])$'s
($V$ an open subgroup of $U$
to which $W$ is sub-conjugate in $G$),
which clearly coincides with the polynomial in \cite{DS2}.
Thus, we complete the proof.
\qed

By definition of $\bar f_U$,
we have
\begin{equation}\label{commutance of res and teich}
\text{\rm Res}_U^G \circ \tau_R=\tau_R \circ \bar f_U.
\end{equation}
Let us define the operator
$\mathcal F_U:R^{\mathcal O(G)}\to R^{\mathcal O(U)}$ by
$$(b_{[V]})_{[V]\in \mathcal O(G)}\mapsto (c_{[W]})_{[W]\in \mathcal O(U)}
$$
where
$$c_{[W]}:=\begin{cases}
b_{[V]} & \text{ if }[W]=[V] \in \mathcal O(G)\,,\\
0 & \text{ otherwise.}
\end{cases}$$
From \cite {O} it follows that
\begin{equation}\label{commutance of fu and Fu}
\mathcal F_U \circ \tilde \varphi
=\tilde \varphi  \circ \text{\rm Res}_U^G .
\end{equation}
In view of Eq. \eqref{commutance of res and teich} and
Eq. \eqref{commutance of fu and Fu}, we have
\begin{equation*}\label{commutance of res and Fu}
\mathcal F_U \circ \Phi
=\Phi \circ  \bar f_U.
\end{equation*}
Consequently, we can complete the following commutative diagrams:
\begin{equation*}\label{restriction0}
\begin{CD}
\mathbb W_G(R) @>{\tau_R}>> \hat \Omega_R(G)\\
@V{\bar f_U}VV @V{\text{\rm Res}_U^G}VV\\
\mathbb W_U(R) @>{\tau_R}>> \hat \Omega_R(U)
\end{CD}
\quad
\begin{CD}
 \hat \Omega_R(G)  @>{\tilde\varphi}>>
R^{\mathcal O(G)} \\
@V{\text{\rm Res}_U^G}VV @V{\mathcal F_U}VV  \\
\hat \Omega_R(U)
@>\tilde \varphi>> R^{\mathcal O(U)}
\end{CD}
\quad
\begin{CD}
\mathbb W_G(R) @>{\Phi}>> R^{\mathcal O(G)} \\
@V{\bar f_U}VV  @V{\mathcal F_U}VV \\
\mathbb W_U(R)
@>\Phi>> R^{\mathcal O(U)}
\end{CD}
\end{equation*}

The above diagrams are also valid with regard to induction operators.
To begin with, by definition of $\bar v_U$,
we have
\begin{equation}\label{commutance of ind and teich}
\text{\rm Ind}_U^G \circ \tau_R=\tau_R \circ \bar v_U.
\end{equation}
It was shown in \cite {O3} that if we define
$$\nu_U:R^{\mathcal O(U)}\to R^{\mathcal O(G)},\quad
(b_{[V]})_{[V]\in \mathcal O(U)} \mapsto (c_{[W]})_{[W]\in \mathcal O(G)},$$
where
\begin{equation*}
c_{[W]}=\sum_{{[V] \in \mathcal O(U)}
\atop {[V]=[W] \text{ in } \mathcal O(G)}}[N_G(W):N_U(V)]\,b_{[V]},
\end{equation*}
then it holds that
\begin{equation}\label{comm of nu and ind}
\tilde \varphi \circ\text{\rm Ind}_U^G=\nu_U\circ \tilde \varphi^{-1}.
\end{equation}
Here, the notation $N_G(W)$ represents the normalizer of $W$ in $G$.
In view of Eq. \eqref{commutance of ind and teich} and
Eq. \eqref{comm of nu and ind}, we have
\begin{equation*}
\nu_U \circ \Phi
=\Phi \circ \bar v_U.
\end{equation*}
Consequently, we have the following commutative diagrams:
\begin{equation*}\label{inductuion0}
\begin{CD}
\mathbb W_U(R) @>{\tau_R}>> \hat \Omega_R(U)\\
@V{\bar v_U}VV @V{\text{\rm Ind}_U^G}VV\\
\mathbb W_G(R) @>{\tau_R}>> \hat \Omega_R(G)
\end{CD}
\quad
\begin{CD}
\hat \Omega_R(U)  @>{\tilde\varphi}>>
R^{\mathcal O(U)} \\
@V{\text{\rm Ind}_U^G}VV @V{\nu_U}VV  \\
\hat \Omega_R(G)
@>\tilde \varphi>> R^{\mathcal O(G)}
\end{CD}
\quad
\begin{CD}
\mathbb W_U(R) @>{\Phi}>> R^{\mathcal O(U)} \\
@V{\bar v_U}VV  @V{\nu_U}VV \\
\mathbb W_G(R)
@>\Phi>> R^{\mathcal O(G)}
\end{CD}
\end{equation*}

\subsection{Necklace polynomials}
In \cite {MR}, Metropolis and Rota introduced the necklace polynomials,
\begin{equation}\label{necklace 0}
\mathcal M(x,n)= \dfrac 1n\displaystyle\sum_{d|n} \mu(d)\,x^{\frac
nd}, \quad n\in \mathbb N,
\end{equation}
and then asked on what class of rings the necklace polynomials in
Eq. \eqref{necklace 0} remain valid. The reason why they
gave this question is because over this class the ring of Witt
vectors becomes isomorphic to the necklace ring.
In this section, we give the answer to this question. To do this, we
will modify the definition of the necklace polynomials in Eq.
\eqref{necklace 0*} and study the properties of the modified ones.

An {\it alphabet} is a set of commuting variables so that,
for example,
$\{x_1,x_2,\cdots,x_m\}$
is the alphabet of variables
$x_1,x_2,\cdots,x_m.$
For alphabets $X=\{x_1,x_2,\cdots,x_m\}$ and
$Y=\{y_1,y_2,\cdots,y_n\}$,
$X \cdot Y$ denotes an alphabet $\{x_ay_b
\,:\, 1\le a \le m, 1\le b \le n\}$.
The notation $\Psi^r(X)$, which is introduced to
be consistent with that of a $\ld$-ring,
will denote an alphabet
$\{x_1^r,x_2^r,\cdots,x_m^r\}$.
The elements of an alphabet $X$ is called {\it letters},
and a {\it word} of $X$ is a finite juxtaposition of letters of $X$.
The {\it length} of a word is the number of letters, where the {\it
product} of two or more words is juxtaposition.
Two words $w$ and $w'$ are said to be {\it conjugate}
when $w=uv$ and $w'=vu$ where $u$ and $v$ are words.
The identity in the monoid of words is the
empty words. An equivalence class of words under the equivalence
relation of conjugacy will be called a {\it necklace}.
If $w=u^i$, then we say that the word of $w$ has {\it period} $n/i$,
where $n$ is the length of $w$.
The smallest $j$ such that $w=v^{n/j}$ for some $v$ is called
the {\it primitive period} of the word $w$.
A word of primitive period $n$ is said to be {\it aperiodic},
and an equivalence class of aperiodic words will be called a {\it
primitive necklace}.

Given a word $w$, let $m(w)$ be the monomial
$x_1^{a_1}x_2^{a_2}\cdots x_m^{a_m}$
where $a_i$ is the number of $x_i$'s in $w$.
Define the
{\it necklace polynomial of degree $k$ over $X$ }
to be
\begin{equation}\label{definition of necklace ring}
M(X, k):=\sum_{w}m(w)\,,
\end{equation}
where the sum is over the primitive necklaces $w$ out of $X$
such that the degree of $m(w)$ equals $k$.
Then, one easily computes
\begin{equation*}
M(X,k)= \dfrac 1k\displaystyle\sum_{d|k} \mu(d)\,p_d(X)^{\frac
kd},
\end{equation*}
where $p_d$ is the $d$-th power sum, that is,
$p_d(X)=x_1^d+\cdots +x_m^{d}.$
For positive integers $i$ and $j$,
we let $[i,j]$ be the least common multiple
and $(i,j)$ the greatest common divisor of $i$ and $j.$
With this notation,
one can get the identities analogous to those in
\cite [Section 3]{MR}.

\begin{thm}\label{necklace 1}\hfill

{\rm (a)} For alphabets $X=\{x_1,x_2,\cdots,x_m\}$ and
$Y=\{y_1,y_2,\cdots,y_n\},$
\begin{align*}
M(X \cdot Y , k)
=\sum_{[i,j]=k}(i,j)
 M(\Psi^{\frac ki}(X) , i) M(\Psi^{\frac kj}(Y), j)\,.
\end{align*}
{\rm (b)} For an alphabet $X=\{x_1,x_2,\cdots,x_m\},$
\begin{align*}
 M(X^r ,k)
=\sum_{j\ge 1 \atop [r,j]=kr}
\dfrac {j}{k}  M(\Psi^{\frac {kr}{j}}(X), j).
\end{align*}
{\rm (c)}
For an alphabet $X=\{x_1,x_2,\cdots,x_m\}$ and
$Y=\{y_1,y_2,\cdots,y_n\}$, we have
\begin{align*}
&(r,s) M(X^{\frac {s}{(r,s)}}Y^{\frac{r}{(r,s)}},k)\\
&=\sum_{i,j} (ri,sj)
M(\Psi^{\frac {[\frac{s}{(r,s)},i]}{i}}(X), i)
M(\Psi^{\frac {[\frac{r}{(r,s)},j]}{j}}(Y), j),
\end{align*}
where $i,j$ range over positive integers such that $ij / (ri,sj)=k / (r,s).$
\end{thm}

\noindent{\bf Proof.}
The proof can be done by a slight modification of
that in \cite[Section 3]{MR}.
So, we will prove only (a).
For a primitive word $w$ of degree $n$ out of the alphabet $XY$,
write it as a monomial in $x_i,y_j$'s, say
$$x_1^{a_1}x_2^{a_2}\cdots x_m^{a_m}y_1^{b_1}y_2^{b_2}\cdots y_n^{b_n}.$$
Let $w'$ be the word $x_1^{a_1}x_2^{a_2}\cdots x_m^{a_m}$ and
$w''$ be the word $y_1^{b_1}y_2^{b_2}\cdots y_n^{b_n}$.
From \cite [Theorem 1, Section 3]{MR} it follows that
the word $w$ is primitive if and only if
$w'$ has primitive period $i$,
the word $w''$ has primitive period $j$, and $[\,i,j\,]=n.$
Write $w'=u'^{\frac ki}$ and $w''=u''^{\frac kj}$ for some $i,j$
satisfying $[i,j]=n.$
Thus, we have $w=u'^{\frac ki}u''^{\frac kj}.$
Now, the bijectivity of this correspondence implies our assertion.
\qed

Another important interpretation of our necklace polynomial
$M(X,k)$ can be described in the following manner.
For a function $f$ from $\mathbb Z /k \mathbb Z$ to $X$,
we let
$$[f]:=\prod_{i=0}^{k-1}f(\bar i).$$
Let us define $\mathbb Z/k\mathbb Z$-action on $f$ by
$\bar  n\cdot f (\bar i):=f(\overline {i-n})$.
Then, the following proposition is almost straightforward.
\begin{prop}
With the above notation, we have
$$ M(X ,k)=\sum_f [f],$$
where the sum is over $f$'s on which $\mathbb Z/k\mathbb Z$
acts freely.
\end{prop}

\begin{rem}
The concept of our necklace polynomial has a nice generalization
associated with a profinite group $G$ and its open subgroups $V$.
In detail, given a datum $(G,V,X)$, the polynomial $M_G(V,X)$ such
that $M_{\hat C}(\hat C^k,X)=M(X,k)$ was introduced. Here, $X$
represents an alphabet. For the complete information refer to
$($\cite{O3}$)$.
\end{rem}

Let $R$ be a special $\ld$-ring. For $r \in R,\, n\in \mathbb N,$
we let
\begin{equation*}
M(r,n):=\dfrac 1n \sum_{d|n}\mu(d)\Psi^{d}(r^{\frac nd}).
\end{equation*}
In order to show this definition to be consistent with
Eq. \eqref{definition of necklace ring},
write $r$ as a sum of one-dimensional elements, say,
$r=r_1+r_2+\cdots+r_m$.
In fact, this expression is possible by virtue of
the {\it splitting principle} for special $\ld$-rings.
Now, consider the alphabet $X_r$ consisting of
$r_i$'s, $ 1\le i \le m $.
Then one can easily verify that
$$M(r,n)=M(X_r,n),$$
and which says that if $R$ is a binomial ring, then
$M(r,n)$ coincides with $\mathcal M(r,n)$
(see \cite{O}).

\begin{thm} Let $R$ be a unital commutative torsion-free ring.
Then, $R$ is a binomial ring if and only if
\begin{equation*}
\mathcal M(r,n)\in R
\end{equation*}
for all $r\in R$ and $n\in \mathbb N$.
\end{thm}

\noindent{\bf Proof.}
The `` if " part follows from Corollary \ref{psiring}.
More precisely, letting $n$ be a prime $p$, the assumption implies that
\begin{equation*}
\mathcal M(r,p)=\frac 1p(r^p-r)\in R.
\end{equation*}
Hence, $R$ satisfies the condition in Corollary \ref{psiring}.
For the `` only if " part, assume that
$r=r_1+r_2+\cdots+r_m$, where $r_i$'s are 1-dimensional.
From the fundamental theorem of symmetric functions it follows that
the symmetric function
\begin{equation*}
M(X,k)= \dfrac 1k\displaystyle\sum_{d|k} \mu(d)\,p_d(X)^{\frac
kd}
\end{equation*}
is an integral polynomial in the elementary symmetric functions
defined using the product of $(1+x_it)$'s.
Here, $X$ denotes the alphabet $\{x_1,x_2,\cdots,x_m\}.$
Now, the specialization that $x_i=r_i, 1\le i\le m$, yields the desired result
since the $n$-th elementary symmetric function is same to
$\ld^n(r)$'s after this specialization.
\qed

\subsection{Covariant functor $Nr$ and its Verschiebung and Frobenius operators }
In this section, we will show that
the functor, $Nr_{\hat C}$, is isomorphic to
the functors $\mathbb W$ and $\Lambda_1$
if these are viewed as the functors from
the category of special $\ld$-rings into itself.
As before, $\hat C$ denotes the profinite completion of
the multiplicative infinite cyclic group $C$.
From now on, we use the notation $Nr$ instead of $Nr_{\hat C}$.
\begin{rem}
(a)
It should be noted that our notation $Nr$ is different from
the one in {\rm \cite{MR}}. Actually, the latter
coincides with the functor $\widehat{Nr}_{\hat C}$ in {\rm \cite{O2}}.

(b)
Very often, the terminologies ``Verschiebung and Frobenius operators"
are used instead of the inductions and the restrictions in case $G=\hat C$.
\end{rem}

We begin with investigating
the structure of the necklace ring $Nr(R)$ over $R$ in more detail,
where $R$ is a special $\ld$-ring.
$Nr(R)$ is the ring whose underlying set is $\prod_{\mathbb N}R$, whose
addition is defined componentwise, and whose
multiplication is given by
\begin{equation}\label{multi of necklace}
(b \cdot c)_n:=\sum_{[i,j]=n}(i,j)\, \Psi^{\frac ni}(b_i)\Psi^{\frac nj}(c_j)
\end{equation}
for two sequences $b=(b_1,b_2,\cdots,b_n,\cdots)$
and $c=(c_1,c_2,\cdots,c_n,\cdots).$
Especially, if $R$ is a binomial ring, then
the multiplication in \eqref{multi of necklace} reduces to
\begin{equation*}
(b \cdot c)_n=\sum_{[i,j]=n}(i,j)\cdot b_i c_j.
\end{equation*}

\begin{lem} \label{isom of neck and groth}
{\rm (\cite [Proposition (17.2.9)]{H})}
For every commutative ring $A$ with identity, the map
\begin{equation}\label{isom E}
E_A:\mathbb W(A) \to \Lambda_1(A), \quad (a_n)_{n\ge 1} \mapsto
\prod_{n=1}^{\infty}\left(\frac{1}{1-a_nt^n}\right)
\end{equation}
is a ring isomorphism.
\end{lem}
For a special $\ld$-ring $R$ let us consider the map
\begin{equation}\label{symmetric power function}
\tilde s_t: Nr(R) \to \Lambda_1(R),\quad
(b_1,b_2,\cdots)\mapsto \prod_{n=1}^{\infty}
\left(\sum_{r=0}^{\infty} S^r(b_n)t^{nr}\right).
\end{equation}
As mentioned in Eq. \eqref{symmetric map}, the notation $S^n$
represents the $n$-th symmetric power operation associated with the
given special $\ld$-ring structure.
It has been shown in \cite {O} that $\tilde s_t$
is a bijective map, and a ring homomorphism if $R$ is
torsion-free.
Indeed, this result holds without the condition ``torsion-free".
\begin{thm}
Let $R$ be a special $\ld$-ring.
Then, $\tilde s_t$ is a ring
isomorphism.
\end{thm}
\noindent{\bf Proof.} It follows from \cite [Lemma 3.13]{O} that
\begin{equation}\label{commuting}
\tilde s_t \circ \tau_R =E_R
\end{equation}
(refer to Eq. \eqref{isom E}).
Combining Theorem \ref{ringhomo of teichmuller}
and Lemma \ref{isom of neck and groth} implies the desired assertion.
\qed
\begin{rem}
Let $R$ be a binomial ring.
Then, the mapping \eqref{symmetric power function} reduces to
$$\tilde s_t (b_1,b_2,\cdots)=\prod_{n=1}^{\infty}
\left(\dfrac {1}{1-t^n}\right)^{b_n},$$ which was dealt with
intensively in $($\cite {MR}$)$.
\end{rem}

We now can make $Nr(R)$ and $\mathbb W(R)$ into special $\ld$-rings
by virtue of the isomorphisms $\tilde s_t$ and $E_A$.
The $\ld^m$-operation on $Nr(R)$ is defined by
$$\ld_m \left({\tilde s_t}^{-1}(f(t))\right)={\tilde s_t}^{-1}(\bar\Lambda_mf(t)).$$
Since
\begin{equation*}\label{recall: def of lambda operation}
\bar \Lambda_m \left(\prod_{i=1}^{\infty}\frac{1}{1-x_it}\right)
=\prod_{i_1<i_2<\cdots <i_m}\frac{1}{1-x_{i_1}\cdots x_{i_m}t}
\end{equation*}
and
\begin{equation}\label{method of commputing lambda}
\tilde s_t(M(x_i))=\frac{1}{1-x_it},
\end{equation}
the $\ld_m$-operation on $Nr(R)$ must satisfy
\begin{equation}\label{lambda operation of necklace}
\ld_m(\sum_{n=1}^{\infty}M(x_i))= \sum_{i_1<i_2<\cdots <i_m}M(x_{i_1}\cdots x_{i_m}).
\end{equation}
Actually, we can verify easily that Eq. \eqref{lambda operation of necklace}
determines the $\ld$-operations completely.
For example, let
$$(c_1,c_2,c_3,\cdots):=\ld^2(b_1,b_2,b_3,\cdots).$$
In order to compute $c_1$
we may assume that $b_1=x_1+x_2$ and $x_i=0$ for $i\ge 3$.
Then, from Eq. \eqref{method of commputing lambda} it follows that
\begin{align*}
b_2=\frac 12(x_1^2-\Psi^2(x_1))+\frac 12(x_2^2-\Psi^2(x_2))
=\frac12 (b_1^2-2c_1-\Psi^{2}(b_1)).
\end{align*}
So, we conclude that $c_1=-b_2+\frac12 (b_1^2-\Psi^{2}(b_1)).$

Similarly, the $\ld_m$-operation on $\mathbb W(R)$ is defined by
$$\ld_m{E_R}^{-1}(f(t))={E_R}^{-1}(\bar\Lambda_mf(t)).$$
Then, from Eq. \eqref{lambda operation of necklace} it follows that
\begin{equation}\label{lambda operation of witt ring}
\ld_m(\sum_{n=1}^{\infty}(x_i,0,0,\cdots))
= \sum_{i_1<i_2<\cdots <i_m}(x_{i_1}\cdots x_{i_m},0,0,\cdots),
\end{equation}
where the summation is being done in $\mathbb W(R).$
Note that Eq. \eqref{lambda operation of witt ring} also
determines the $\ld$-operations of $\mathbb W(R)$ completely.
For example, letting
$$(d_1,d_2,d_3,\cdots):=\ld^2(e_1,e_2,e_3,\cdots),$$
then it is easy to show that $d_1=-e_2.$

Now, let us investigate morphisms.
For a special $\ld$-ring homomorphism,
$f:A \to B$, consider the homomorphisms
\begin{align*}
&\Lambda(f):\Lambda_1(A)\to \Lambda_1(B),\quad
1+\sum a_nt^n\mapsto 1+\sum f(a_n) t^n,\\
&Nr(f): Nr(A)\to Nr(B),\quad (b_n)_n \mapsto (f(b_n))_n,\\
&\mathbb W(f): \mathbb W(A)\to \mathbb W(B),\quad (a_n)_n \mapsto (f(a_n))_n.
\end{align*}
It is well-known that $\Lambda(f)$ is a special $\ld$-ring homomorphism.
Therefore, $Nr(f)$ and $\mathbb W(f)$ also are special $\ld$-ring homomorphisms
by the definition of their $\ld$-operations.
In addition, following the same way as in \cite{O2}, we can show that
\begin{align*}
&\tau_B \circ \mathbb W(f)\circ \tau_A^{-1}
=\tilde s_t^{-1} \circ \Lambda(f)\circ \tilde s_t,\\
&Nr(f)=\tau_B \circ \mathbb W(f)\circ \tau_A^{-1}.
\end{align*}
The discussion until now can be illustrated in
the following commutative diagram:
\begin{equation*}\label{SMALL diagram}
\begin{aligned}
\begin{array}{ccccc}
\mathbb W(A)
& \stackrel{\tau_A}{\longrightarrow}
&Nr(A)
& \stackrel{\tilde s_{t}}{\longrightarrow}
&\Lambda_1(A)\\
\Big\downarrow\vcenter{\rlap{$\scriptstyle{\mathbb W(f)}$}}
&{}
&\Big\downarrow\vcenter{\rlap{$\scriptstyle{Nr(f)}$}}
&{}
&\Big\downarrow\vcenter{\rlap{$\scriptstyle{\Lambda(f)}$}}\\
\mathbb W(B)
& \stackrel{ \tau_{B}}{\longrightarrow}
&Nr(B)
& \stackrel{ \tilde s_{t}}{\longrightarrow}
&\Lambda_1(B)
\end{array}
\end{aligned}
\end{equation*}
Next, let us discuss Verschiebung and Frobenius operators on
$Nr(R)$ and $\Lambda_1(R)$ more closely.
Indeed, they are nothing but inductions and restrictions respectively
with $G=\hat C$.
The $r$-th Verschiebung operator, $V_r$, on $Nr(R)$
is defined to be
$$V_r \alpha=\beta,\text{ where }
\beta_n=\alpha_{\frac nr}
\text{ with } \alpha_{\frac nr}:=0 \text{ if } n/r \notin \mathbb N.$$
The $r$-th Frobenius operator, $F_r$, on $Nr(R)$ is defined to be
$$F_r \alpha=\beta,\text{ where }
\beta_n=\sum_{[r,j]=rn}\frac jn  \Psi^{\frac {nr}{j}}(\alpha_j).$$
Note that $F_r$ is a ring isomorphism, whereas $V_r$ is just additive.
For $r\in R$, we let
\begin{equation*}
M(r):=(M(r,1), M(r,2),\cdots,  M(r,n),\cdots)
\end{equation*}
(refer to \cite{O}).
\begin{prop} {\rm (cf. \cite {MR})}\label{properties of ver and Frob}
For $a,b\in R$, $\alpha\in Nr(R)$, and $r,s\in \mathbb N$, we have

{\rm (a)}
$V_r V_r=V_{rs}.$\\
{\rm (b)}
$F_r F_s=F_{rs}.$\\
{\rm (c)}
$F_r V_r(\alpha)=r\alpha.$\\
{\rm (d)}
$V_r M(a) \cdot  V_s  M(b)=(r,s) V_{[r,s]}
M(a^{\frac {r}{(r,s)}}b^{\frac {s}{(r,s)}})\,.$\\
{\rm (e)} $F_r M(a)= M(a^r)$.\\
{\rm (f)} $F_r V_s=(r,s)V_{\frac {[r,s]}{r}}F_{\frac {[r,s]}{s}}$.
\end{prop}

\noindent{\bf Proof.}
(a) and (b) follow from Eq. \eqref{transitivity of ins and res}.

(c) Letting $b=F_rV_r(\alpha)$, then
\begin{equation*}
b_n
=\sum_{[r,j]=nr}\dfrac jn \Psi^{\frac {nr}{j}}(V_r(\alpha)_j).
\end{equation*}
For $j=ri$ the condition $[r,j]=nr$ implies that $i=n$.
Thus, we have the desired result.

(d)
Letting $c=(c_n)_n=V_r M(\alpha) \cdot V_s M(\beta),$ then
\begin{align*}
c_n&=\sum_{[i,j]=n}(i,j)
\Psi^{\frac ni}(M(\alpha; \frac ir))\Psi^{\frac nj}(M(\beta;\frac js)).
\end{align*}
If we substitute $\frac ir=i'$ and $\frac js=j'$,
then the desired result is immediate.

(e) By definition of $F_r$ we have
\begin{align*}
F_r M(a)&=\sum_{[r,j]=rn}\frac jn \Psi^{\frac {nr}{j}}(M(a;j))\\
&=M(a^r,j) \quad \text{\rm by Theorem \ref{necklace 1}(b).}
\end{align*}

(f) This identity follows from \cite[Proposition 3.10]{O}.
\qed

We now suppose that $R$ has two different kinds of special $\ld$-ring structures.
Write their $\ld$-operations and Adams operations as
$$(\ld_1^n, \Psi_1^n ; n\ge 1)\quad
\text { and } \quad
(\ld_2^n, \Psi_2^n ; n\ge 1).$$
In Corollary \ref{canonical isom}
we have shown that there exists a canonical isomorphism
\begin{align*}
& Nr({\rm id}_{12})=\tau_{R,2} \circ{\tau_{R,1}}^{-1}:Nr_1(R) \to Nr_2(R),\\
&\sum_{n=1}^{\infty} V_n M_1(q_n) \mapsto \sum_{n=1}^{\infty} V_n M_2(q_n)\,,
\end{align*}
where
$$M_i(r):=(M_i(r,1),M_i(r,2),\cdots)$$
with
$$M_i(r,k)=\dfrac 1k \sum_{d|k}\mu(d)\Psi_i^{d}(r^{\frac kd}).$$
The map $Nr({\rm id}_{12})$
behaves nicely for Verschiebung and Frobenius operators.
\begin{prop}\label{operator-preserving}
The isomorphism $Nr({\rm id}_{12})$ preserves Verschiebung and Frobenius operators,
that is,

{\rm (a)}
$Nr({\rm id}_{12})(V_r \,\alpha)
=V_r Nr({\rm id}_{12})(\alpha),\quad  \alpha\in Nr_1(R).$\\
{\rm (b)}
$Nr({\rm id}_{12})(F_r \,\alpha)
=F_r Nr({\rm id}_{12})(\alpha),\quad \alpha\in Nr_1(R).$
\end{prop}
\noindent{\bf Proof.}
Since $\tau_{R,i},\,i=1,2$ preserve
inductions and restrictions
by Eq. \eqref{commutance of res and teich} and
Eq. \eqref{commutance of ind and teich},
so does $Nr({\rm id}_{12})$.
\qed

The $n$-th Verschiebung operator, $V^\Lambda_n$, is defined by
\begin{equation*}\label{ver of Grothen}
V^\Lambda_n(1+a_1t+a_2t^2+\cdots):=1+a_1t^n+a_2t^{2n}+\cdots.
\end{equation*}
In order to define the $n$-th Frobenius operator, $F^\Lambda_n$,
we define $a_i$'s
by the equation
$$1+a_1 t +a_2t^2+\cdots=\prod_{i=1}^{\infty}\frac {1}{1-x_i t}.$$
Let $Q_{n,k}(a_1,\cdots,a_{nk})$ be the coefficient of $t^k$ in
$$\prod_{i=1}^{\infty}\frac{1}{1-x_i^nt}.$$
Now, we define
\begin{equation*}\label{frob of Grothen}
F_n^\Lambda(1+\sum c_it^i)=1+\sum_{k=1}^{\infty} Q_{n,k}(c_1,\cdots,c_{nk})t^k.
\end{equation*}

\begin{prop}{\rm (cf. Hazewinkel \cite {H})}\label{Frob=res}
The $n$-th Frobenius operator $F_n^\Lambda$ coincides with
the $n$-th Adams operator $\Psi^{n}$ .
\end{prop}
\noindent{\bf Proof.}
This assertion was proved for the ring $\Lambda_0(A)$ in \cite {H}.
Hence, applying the isomorphism
\begin{equation*}\label{isom of grothendict rings}
\iota : \Lambda_0(A)\stackrel{\cong}{\to} \Lambda_1(A),
\quad f(t)\mapsto \frac {1}{f(t)}.
\end{equation*}
we obtain the desired result.
\qed

\begin{prop}
Let $R$ be a special $\ld$-ring.
Then, the map $\tilde s_t$ in Eq. \eqref{symmetric power function}
preserves Verschiebung and Frobenius operators.
\end{prop}
\noindent{\bf Proof.}
In order to prove the assertion, it is enough to show that
\begin{align*}
&\tilde s_t(V_r V_s M(\alpha))=V_r^\Lambda \tilde s_t(V_s M(\alpha)),\\
&\tilde s_t(F_r V_s M(\alpha))=F_r^\Lambda \tilde s_t(V_s M(\alpha)).
\end{align*}
For the first identity, let us combine Eq. \eqref{commuting}
with Proposition \ref{properties of ver and Frob} (a) to get
\begin{align*}
\tilde s_t(V_r V_s M(\alpha))=\tilde s_t(V_{rs} M(\alpha))
=\dfrac {1}{1-\alpha t^{rs}}=V_r^\Lambda \dfrac{1}{1-\alpha t^s}.
\end{align*}
Since $\tilde s_t(V_sM(\alpha))=\dfrac{1}{1-\alpha t^s}$, we are done.
For the second one, note that
\begin{align*}
\tilde s_t(F_r\circ V_s M(\alpha))
&=\tilde s_t((r,s)V_{\frac {[r,s]}{r}}F_{\frac {[r,s]}{s}}M(\alpha))\\
&=\tilde s_t((r,s)V_{\frac {[r,s]}{r}}M(\alpha^{\frac {[r,s]}{s}}))
\quad \text{ (by Proposition \ref{properties of ver and Frob} (e))}\\
&=\left(\dfrac {1}{1-\alpha^{\frac {[r,s]}{s}}t^{\frac {[r,s]}{r}}}\right)^{(r,s)}.
\end{align*}
On the other hand,
\begin{align*}
F_r^\Lambda \tilde s_t(V_s M(\alpha))
&=F_r^\Lambda V_s^\Lambda \tilde s_t(M(\alpha))\\
&=(r,s)V^\Lambda_{\frac {[r,s]}{r}}F^\Lambda_{\frac {[r,s]}{s}}\tilde s_t(M(\alpha))\\
&=(r,s)V^\Lambda_{\frac {[r,s]}{r}}\left(\dfrac {1}{1-\alpha^{\frac {[r,s]}{s}}t}\right)\\
&=\left(\dfrac {1}{1-\alpha^{\frac {[r,s]}{s}}t^{\frac {[r,s]}{r}}}\right)^{(r,s)}.
\end{align*}
This completes the proof.
\qed

The results in this section may be summarized as follows.
The diagram
\begin{equation}\label{big diagram}
\begin{aligned}
\begin{array}{ccccc}
\mathbb W(R)
& \stackrel{\tau_{R}(\cong)}{\longrightarrow}
&Nr(R)
& \stackrel{\tilde s_t(\cong)}{\longrightarrow}
&\Lambda_1(R)\\
\Big\downarrow\vcenter{\rlap{$\scriptstyle{\Phi}$}}
&{}
&\Big\downarrow\vcenter{\rlap{$\scriptstyle{\tilde \varphi}$}}
&{}
&\Big\downarrow\vcenter{\rlap{$\scriptstyle{\frac {d}{dt}\log}$}}\\
{\rm gh}(R)
&\stackrel{\text{id}}{\longrightarrow}
&{\rm gh}(R)
&\stackrel{\text{\rm identification}}{\longrightarrow}
&R[[t]]
\end{array}
\end{aligned}
\end{equation}
is commutative and all the maps appearing in this diagram
preserve Verschiebung and Frobenius operators.
Here, ${\rm gh}(R)$ represents the ghost ring of $R$, which is the set
$\prod_{\mathbb N}R$ with the addition and the multiplication
defined componentwise.

We end this section by introducing two significant
properties of $\tilde \varphi$.

The first is that
for all $n\ge 1$ the maps
$\tilde \varphi_n$ provide natural transformations from the functor
$Nr$ to the identity functor, which follows
from the commutative diagram
\begin{equation*}
\begin{aligned}
\begin{array}{ccc}
Nr(A)
& \stackrel{Nr(f)}{\longrightarrow}
&Nr(B)\\
\Big\uparrow\vcenter{\rlap{$\scriptstyle{\tau_A}$}}
&{}
&\Big\uparrow\vcenter{\rlap{$\scriptstyle{\tau_B}$}}
\\
\mathbb W(A)
&\stackrel{\mathbb W(f)}{\longrightarrow}
&\mathbb W(B)\\
\Big\downarrow\vcenter{\rlap{$\scriptstyle{\Phi_n}$}}
&{}
&\Big\downarrow\vcenter{\rlap{$\scriptstyle{\Phi_n}$}}
\\
A& \stackrel{f}{\longrightarrow}
&B
\end{array}
\end{aligned}
\end{equation*}
Here, $\Phi_n$ is the projection of $\Phi$ to the $n$-th component.
Observe that
\begin{align*}
f\circ \tilde \varphi_n
&=f\circ \Phi_n \circ \tau_A^{-1}\\
&=\Phi_n \circ \mathbb W(f)\circ \tau_A^{-1}\\
&=\Phi_n \circ \tau_B^{-1} \circ Nr(f)\\
&=\tilde \varphi_n \circ Nr(f)\,.
\end{align*}

The second is related to a generalization
of Theorem \ref{necklace 1} {\rm (a)}.
Let $R$ be a $\mathbb Q$-algebra and
${\bf a}=(a_n)_n,\, {\bf b}=(b_n)_n \in \prod_{\mathbb N}R$.
Set
\begin{equation*}
E({\bf a},n)=\dfrac 1n \sum_{d|n}\mu(d)\Psi^{d}(a_{\frac nd}).
\end{equation*}
Since $\tilde \varphi$ is a ring homomorphism
we can derive the identity
\begin{equation*}
E({\bf ab}, n)
=\sum_{[i,j]=n}(i,j)
\Psi^{\frac ni}(E({\bf a},i))\Psi^{\frac nj}(E({\bf b},j)).
\end{equation*}
In particular, considering the case
${\bf a}=(r,r^2,r^3,\cdots)$ and ${\bf b}=(s,s^2,s^3,\cdots)$,
we can recover Theorem \ref{necklace 1} {\rm (a)}.

\section
{LOGARITHMIC FUNCTIONS ASSOCIATED WITH GRADED LIE SUPERALGEBRAS}
\label{LOGARITHM}
\subsection
{Definition}\label{Definition:Logarithm}

Let $\widehat {\Gamma}$ be a free abelian group with finite rank
and let $\Gamma$ be a countable (usually infinite) sub-semigroup
in $\widehat {\Gamma}$ such that every element $\alpha \in \Gamma$
can be written as a sum of elements of $\Gamma$ in only finitely
many ways.
Given a commutative algebra $R$ with unity over $\mathbb C$,
consider
$$R [[\G \times \mathbb Z_2]]= \{\sum_{(\la,a) \in \G\times
\mathbb Z_2} \zeta(\la,a) E^{(\la,a)}\, : \, \zeta(\la,a) \in R
\},$$
the completion of the semi-group algebra
$R[\G \times \mathbb Z_2]$.
Here $E^{(\la,a)}=(-1)^a e^{(\la,a)}$ and
$e^{(\la,a)}$ are the usual basis elements of
$R[\G\times \mathbb Z_2]$ with the multiplication given by
$e^{(\lambda,a)}e^{(\mu,b)}=e^{(\lambda+\mu,a+b)}$ for
$(\lambda,a),(\mu,b) \in \G\times \mathbb Z_2$.
Then it is easy to show that
$\{E^{(\la,a)}: (\la,a) \in \G \times \mathbb Z_2\}$ is also a basis
of $R[\G\times \mathbb Z_2]$ with the multiplication
$E^{(\lambda,a)}E^{(\mu,b)}=E^{(\lambda+\mu,a+b)}.$

For each $\alpha=\sum_{i=1}^{r}k_i\alpha_i \in \G,$ we define the
{\it height} of $\alpha$, denoted by ${\rm ht}(\alpha)$, to be the
number $\sum_{i=1}^r k_i.$ We write $R [[\G\times \mathbb Z_2]]^{(n)}$
for the set consisting of all elements
$$\sum\zeta(\beta,b)E^{(\beta,b)},
\text { where }\zeta(\beta,b)=0 \text{ whenever } {\rm ht}(\beta) <n.$$
It is clear that $R [[\G\times \mathbb Z_2]]^{(n)}$ is an ideal of $R
[[\G\times \mathbb Z_2]],$ and therefore we have a filtration of
ideals
$$R [[\G\times \mathbb Z_2]]=  R [[\G\times \mathbb Z_2]]^{(1)}
\supset R [[\G\times \mathbb Z_2]]^{(2)}\supset R [[\G\times \mathbb
Z_2]]^{(3)}\supset\cdots.$$
For $f\in R [[\G\times \mathbb Z_2]]$ and
$c\in \C$, define $(1+f)^c$ by
$$1+cf+\frac{c(c-1)}{2!}f^2+\frac{c(c-1)(c-2)}{3!}f^3+\cdots.$$

\begin{df} {\rm(cf. \cite {B})} We call a function
$ L : R [[\G\times \mathbb Z_2]] \to R [[\G\times \mathbb Z_2]]$
{\it logarithmic} if it satisfies the following properties.
\begin{align*}
&L ({R [[\G\times \mathbb Z_2]]}^{(n)}) \subseteq
{R [[\G\times \mathbb Z_2]]}^{(n)}
\text{ for all } n \in \mathbb N,\\
& L (f)+L(g)=L (f+g-fg)
\text { for all } f,g \in R [[\G\times \mathbb Z_2]],\\
&c  L (f)=L (1-(1-f)^c) \text { for all }c\in \C, f\in R
[[\G\times \mathbb Z_2]].
\end{align*}
\end{df}

Typical examples of logarithmic functions are Log and Exp, which
are defined to be
\begin{align}
& \label {log and exp1 }\text{\rm Log}(f)=-\log(1-f)=\sum_{k \ge 1 }\dfrac 1k f^k,\\
& \label {log and exp2}\text{\rm Exp}(f)=1-\exp(-f)=-\sum_{k \ge
1 }\dfrac {(-1)^k}{k} f^k
\end{align}
for all $f\in R [[\G\times \mathbb Z_2]].$
One can easily show that they are mutually inverses to each other.

Given an element $f=\sum \zeta(\alpha,a)E^{(\alpha,a)}\in R
[[\G\times \mathbb Z_2]],$ the coefficient of $E^{(\alpha,a)}$ in
Log$(f)$ can be computed as follows. Let
$$P(f)=\{(\alpha,a)\in \G\times \mathbb Z_2 \,:\,
\zeta(\alpha,a)\ne 0\},$$
and let
$\{(\beta_i,b_j)|\,i,j=1,2,3,\cdots\}$
be a fixed enumeration of $P(f).$
For $(\alpha,a)\in P(f),$ set
$$T(\alpha,a)=\{s=(s_{ij})\,:\,s_{ij} \ge 0, \sum s_{ij}(\beta_i,b_j)=(\alpha,a)\}.$$
And then we set
\begin{equation}\label{coef of logarithmic}
W(\alpha,a)=\sum_{s\in T(\alpha,a)}\dfrac{(|s|-1)!}{s!}
\prod {\zeta(\beta_i,b_j)}^{s_{ij}},
\end{equation}
where
$|s|=\sum_{i,j}s_{ij}$ and $s!=\prod_{i,j}s_{ij}!.$
By direct computation of the right hand side of Eq. \eqref {log and exp1 }
we can derive the following expansion formula:
$$\text{Log}(f)=\sum W(\alpha,a)E^{(\alpha,a)}.$$

Two basic properties of logarithmic functions were given by Bryant
in the case where the rank of $\Gamma$ is one
(see \cite[Theorem 2.2 and 2.3]{B}).
Actually one can easily generalize those properties to the general case
following Bryant's approach.
Let $\{v_{\lambda}|\, \lambda \in \Lambda \}$ be a $\C$-basis of
$R.$ Then, as a $\C$-vector space, $R [[\G\times \mathbb Z_2]]$
has a basis
$$\{v_{\lambda}E^{(\alpha,a)}\,:\, \lambda \in \Lambda, (\alpha,a) \in
\G\times \mathbb Z_2\}.$$
To each $v_{\lambda}E^{(\alpha,a)}$ we assign an arbitrary element
$f_{\lambda,(\alpha,a)}\in R
[[\G\times \mathbb Z_2]]^{({\rm ht} (\alpha))}.$

\begin{prop}\label{basis theorem}\hfill

{\rm (a)} There exists a unique logarithmic function $L: R
[[\G\times \mathbb Z_2]] \to R [[\G\times \mathbb Z_2]]$ such that
$$L(v_{\lambda}E^{(\alpha,a)})=f_{\lambda,(\alpha,a)}
,\quad (\alpha,a)\in \G\times \mathbb Z_2.$$

{\rm (b)} A function $L$ is logarithmic if and only if $L=\Phi \circ
\text{\rm Log}$ for some $\C$-linear function $\Phi$ satisfying
$\Phi(R[[\G\times \mathbb Z_2]]^{(n)})
\subseteq R [[\G\times \mathbb Z_2]]^{(n)}$
for all $n \ge 1.$
\end{prop}

For later use we introduce a shift operator on
$R [[\G\times \mathbb Z_2]].$
Given each positive integer $k,$ let
$$\Theta_k :R[[\G\times \mathbb Z_2]] \to R [[\G\times \mathbb Z_2]],
\quad \sum\zeta(\beta,b)E^{(\beta,b)}
\mapsto \sum\zeta(\beta,b)E^{k(\beta,b)}.$$

Very often we use $\Gamma$-grading instead of $(\G\times
\mathbb Z_2)$-grading.
In this case we let
$$\zeta(\alpha)=\zeta(\alpha, \overline 0) + \zeta(\alpha,
\overline 1),\quad E^{(\alpha,\overline 0)}=E^{(\alpha,\overline 1
)}=E^{\alpha} , \quad (\alpha \in \Gamma).$$
Also, $R [[\G\times \mathbb Z_2]]$ will be replaced by
$$R [[\G]]= \{\sum_{\la \in \G} \zeta(\la) E^{\la} \,: \, \zeta(\la) \in R \}.$$

\subsection{Graded Lie superalgebra and its Lie module denominator identity}

Let $G$ be a group and $K$ be a field with char$K\ne 2.$
Throughout this section we assume that $\G_K(G)$ is a special
$\ld$-ring. Consider a ($\G \times \mathbb Z_2$)-graded $K$-vector
space $V=\bigoplus_{(\alpha,a)\in \G \times \mathbb
Z_2}V_{(\alpha,a)}$ with $\dim V_{(\alpha,a)} < \infty$ for all
$(\alpha,a) \in (\G \times\mathbb Z_2)$. Furthermore, if $G$ acts
on $V$ preserving the $\G\times \mathbb Z_2$-gradation, we define
$$[V]:=
\sum_{(\alpha,a)} [ V_{(\alpha,a)}]e^{(\alpha,a)} \in \G_K(G)[[\G
\times \mathbb Z_2]].$$ If we set
$|V_{(\alpha,a)}|:=(-1)^a [V_{(\alpha,a)}],$ then it is easy to show that
$$[V]
=\sum_{(\alpha,a)} |V_{(\alpha,a)}|E^{(\alpha,a)}.$$ We call $[V]$
the {\it $G$-module function of $V$ over $K.$}

On the other hand, a $\Z_2$-graded $K$-vector space
$\fL=\fL_{\overline 0} \oplus \fL_{\overline1}$ is called a {\it
Lie superalgebra} if there exists a bilinear map [\ ,\ ] : $\fL
\times \fL \to \fL$, called the {\it bracket}, such that
\begin{align*}
& [\fL_a,\fL_b]\subset
\fL_{a+b},\\
& [x,y]=-(-1)^{ab} [y,x],\\
& [x,[y,z]]=[[x,y],z]+(-1)^{ab}[y,[x,z]]
\end{align*}
for all $x \in \fL_{a}, \ y \in \fL_{b}, \ a, b \in \Z_2$. If char
$K=3,$ then additionally
$[x,[x,x]]=0$ for $x \in \fL_{\overline 1}.$

The homogeneous elements of ${\fL}_{\overline 0}$ (resp.
${\fL}_{\overline 1}$) are called {\it even} (resp. {\it odd}).
Consider a ($\G \times \mathbb Z_2$)-graded Lie superalgebra
$\mathfrak L=\bigoplus_{(\alpha,a)\in \G \times \mathbb Z_2}\mathfrak
L_{(\alpha,a)}$ with $\dim \mathfrak L_{(\alpha,a)} < \infty$ for all
$(\alpha,a) \in (\G \times\mathbb Z_2).$
In addition we suppose that $G$
acts on $\mathfrak L$ preserving the ($\G \times \mathbb
Z_2$)-gradation.

\begin{rem}
When we deal with Lie algebras instead of
Lie superalgebras, $K$ may have characteristic $2.$ In this case,
we assume that $[x,x]=0$ for all $x \in \mathfrak L.$ All the
results for Lie superalgebras appearing in this section may be
carried over to Lie algebras over an arbitrary field, in
particular, of characteristic 2.
\end{rem}

For each ($\G \times \mathbb Z_2$)-graded Lie superalgebra
$\mathfrak L,$ the homology groups of $\mathfrak L$ are defined as
the torsion groups of its universal enveloping algebra viewed as a
supplemented algebra (see \cite {KK98}).

Let ${\mathfrak L}={\mathfrak L}_{\overline 0}\oplus {\mathfrak L}_{\overline
1}$ be a Lie superalgebra and $U=U({\mathfrak L})$ be its universal
enveloping algebra. For each $k\geq 0$, define
\begin{equation*} \label {3.1}
C_{k}=C_k({\mathfrak L})=\bigoplus_{p+q=k}\Lambda^p({\mathfrak
L}_{\overline 0}) \otimes S^q({\mathfrak L}_{\overline 1}).
\end{equation*}
Consider the following chain complex
$(M_{k}, \pr_k)$ $(k\geq -1)$, where
$$M_{k} = \begin{cases} U \otimes_K C_{k} ({\mathfrak L}) &
\text {if} \ k\ge 0, \\
K & \text {if} \ k=-1,
\end{cases}$$
and the differentials $\pr_k:M_k \rightarrow M_{k-1}$ are given by
\begin{equation*}
\begin{aligned}
& \pr_{k}(u\otimes(x_1\wed\dotsb\wed x_p)\otimes(y_1\dotsb y_q)) \\
&=\sum_{1\leq s <t \leq p}(-1)^{s+t}
  u\otimes([x_s,x_t]\wed x_1\wed\dotsb\wed\widehat{x_s}\wed\dotsb\wed
\widehat{x_t}\wed\dotsb\wed x_p)
 \otimes(y_1\dotsb y_q) \\
& +\sum_{s=1}^{p}\sum_{t=1}^{q}(-1)^{s}
    u\otimes(x_1\wed\dotsb\wed\widehat{x_s}\wed\dotsb\wed x_p)
    \otimes([x_s,y_t]y_1\dotsb\widehat{y_t}\dotsb y_q)  \\
&-\sum_{1\leq s < t\leq q }
    u\otimes([y_s,y_t]\wed x_1\wed\dotsb\wed x_p)
     \otimes(y_1\dotsb\widehat{y_s}\dotsb\widehat{y_t}\dotsb y_q) \\
&+\sum_{s=1}^{p}(-1)^{s+1}
    (u\cdot x_s)\otimes(x_1\wed\dotsb\wed\widehat{x_s}\wed\dotsb\wed x_p)
    \otimes(y_1\dotsb y_q)  \\
& + (-1)^p\sum_{t=1}^{q}(u\cdot y_t)
   \otimes(x_1\wed\dotsb\wed x_p)\otimes(y_1\dotsb\widehat{y_t}\dotsb y_q)
\end{aligned}
\end{equation*}
for $k\ge 1$, $\pr_{0}$ is the augmentation map extracting the
constant term, and $\pr_{-1}=0$. Then one can easily verify that
$\pr_{k-1} \circ \pr_{k}=0$.

In particular, when ${\mathfrak L}$ is a $\G$-graded Lie algebra and
$U=U({\mathfrak L})$ is its universal enveloping algebra over an
arbitrary field $K$, then
\begin{equation*}
C_{k}=C_k({\mathfrak L})=\Lambda^k({\mathfrak L}),
\end{equation*}
$$M_{k} =
\begin{cases}
U \otimes_K C_{k}({\mathfrak L})\ \
&\text {if} \ k\ge 0, \\
K \ \ & \text {if} \ k=-1,
\end{cases}$$
and the differentials
$\pr_k:M_k \rightarrow M_{k-1}$ reduce to
\begin{equation*}\label{differentials 1}
\begin{aligned}
& \pr_{k}(u\otimes(x_1\wed\dotsb\wed x_k)) \\
&=\sum_{1\leq s <t \leq k}(-1)^{s+t}
  u\otimes([x_s,x_t]\wed x_1\wed\dotsb\wed\widehat{x_s}\wed\dotsb\wed
\widehat{x_t}\wed\dotsb\wed x_k) \\
&+\sum_{s=1}^{k}(-1)^{s+1}
(u\cdot x_s)\otimes(x_1\wed\dotsb\wed\widehat{x_s}\wed\dotsb\wed x_k).
\end{aligned}
\end{equation*}
\begin{prop} {\rm (\cite{CE,KK98})}
\label{Prop 3.1} \label {resolution}\hfill

{\rm (a)} Let $\mathfrak L$ be a Lie algebra over an arbitrary
field $K$. Then, the chain complex $M=(M_{k}, \pr_{k})$ is a free
resolution of the trivial 1-dimensional left $U$-module $K.$

{\rm (b)} Let $\mathfrak L$ be a Lie superalgebra over $\C.$
Then, the chain complex $M=(M_{k}, \pr_{k})$ is a free resolution of the
trivial 1-dimensional left $U$-module $\C.$
\end{prop}

\begin{rem}
In Proposition \ref {resolution} {\rm (b)},
$\C$ may be replaced by a field of characteristic zero.
However, it does not seem to be known yet whether
this statement remains true
over an arbitrary field.
\end{rem}

Let $K$ be the trivial one dimensional ${\mathfrak L}$-module and
consider the chain complex $(K\otimes_U M_k, 1_K \otimes_U
\partial_k ).$
Then, it is straightforward that
the homology modules of this chain complex, denoted by
$H_{k}({\mathfrak L})=H_{k}({\mathfrak L},K),$ are determined from the complex
\begin{equation*}\label {2.8}
\dotsb\rightarrow C_{k}({\mathfrak L}) \stackrel{d_k}{\rightarrow}
C_{k-1}({\mathfrak L})\stackrel{d_{k-1}}\rightarrow\dotsb \rightarrow
C_{1}({\mathfrak L})\stackrel{d_1}{\rightarrow} C_{0}({\mathfrak L})
\rightarrow 0,
\end{equation*}
where $C_{k}({\mathfrak L})$ are defined by
\begin{equation*}
C_{k}({\mathfrak L})=\bigoplus_{p+q=k} \Lambda^{p}({\mathfrak
L}_{\overline 0})\otimes S^{q}({\mathfrak L}_{\overline 1})
\end{equation*}
and the differentials $d_k:C_{k}({\mathfrak L})\rightarrow
C_{k-1}({\mathfrak L})$ are given by
\begin{align*}
& d_{k}((x_1 \wed \dotsb \wed x_p)\otimes(y_1 \dotsb y_q)) \\
& =\sum_{1\leq s < t\leq p}(-1)^{s+t} ([x_s,x_t]\wed
x_1\wed\dotsb\wed \widehat{x_s}\wed
\dotsb\wed\widehat{x_t}\wed\dotsb\wed x_p)
\otimes(y_1\dotsb y_q) \\
&+\sum_{s=1}^{p}\sum_{t=1}^{q}(-1)^{s}
(x_1\wed\dotsb\wed\widehat{x_s}\wed\dotsb\wed x_p)\otimes
([x_s,y_t]y_1\dotsb\widehat{y_t}\dotsb y_q) \\
 &-\sum_{1\leq s < t\leq q}
([y_s,y_t]\wed x_1\wed\dotsb\wed x_p)
\otimes(y_1\dotsb\widehat{y_s}\dotsb\widehat{y_t}\dotsb y_q)
\end{align*}
for $k\geq 2$, $x_i\in {\mathfrak L}_{\overline 0}$, $y_j\in {\mathfrak
L}_{\overline 1}$, and $d_1=0.$

Let
\begin{align*}
& C({\mathfrak L})=\sum_{k=0}^{\infty}(-1)^{k}C_{k}({\mathfrak L})
=\mathbb{C}\ominus \mathfrak L\oplus C_{2}({\mathfrak L})\ominus\dotsb , \\
& \Lambda({\mathfrak L}_{\overline 0})=\sum_{k=0}^{\infty}(-1)^{k}
\Lambda^{k}({\mathfrak L}_{\overline 0}) =\mathbb{C}\ominus {\mathfrak
L}_{\overline 0} \oplus
\Lambda^{2}({\mathfrak L}_{\overline 0})\ominus\dotsb , \\
& S({\mathfrak L}_{\overline
1})=\sum_{k=0}^{\infty}(-1)^{k}S^{k}({\mathfrak L}_{\overline 1})
=\mathbb{C}\ominus {\mathfrak L}_{\overline 1} \oplus S^{2}({\mathfrak
L}_{\overline 1})
\ominus\dotsb , \\
& H({\mathfrak L})=\sum_{k=1}^{\infty}(-1)^{k+1}H_{k}({\mathfrak L})
=H_1({\mathfrak L})\ominus H_2({\mathfrak L})\oplus H_3({\mathfrak
L})\ominus\dotsb
\end{align*}
be the alternating direct sums of $(\G\times \mathbb Z_2)$-graded
vector spaces. Clearly
$$C({\mathfrak L})=\Lambda({\mathfrak L}_{\overline 0})
\otimes S({\mathfrak L}_{\overline 1}).$$

Let $t=(t\ag)_{\ag \in \GG}$ be a sequence of nonnegative integers
indexed by ($\GG$) with only finitely many nonzero terms, and set
$|t|=\sum t\ag$. Since the $k$-th exterior power $\Lambda^k
({\mathfrak L}_{\overline 0})$ is decomposed as
\begin{equation*}
\Lambda^k({\mathfrak L}_{\overline 0})=
\bigoplus_{|t|=k}\left(\bigotimes_{\al \in \G}
\Lambda^{t(\al,{\overline 0})}({\mathfrak L}_{(\al,{\overline
0})})\right)
\end{equation*}
as a $(\G\times \mathbb Z_2)$-graded $KG$-module, we have
\begin{align*}
[\Lambda({\mathfrak L}_{\overline 0})]
&=\sum_{k=0}^{\infty}(-1)^k [\Lambda^k({\mathfrak L}_{\overline 0})] \\
&=\prod_{\al \in \G} \left(\sum_{m=0}^{\infty}(-1)^m
[\Lambda^m({\mathfrak L}_{(\al,\overline 0})]e^{m(\al,{\overline
0})}\right).
\end{align*}
By the definition of Adams operations, we have
\begin{equation*}
 [\Lambda({\mathfrak L}_{\overline 0})]=\prod_{\al \in \G} \exp
\left(-\sum_{r=1}^{\infty}\frac{1}{r}\Psi^r([{\mathfrak
L}_{(\al,{\overline 0})}])e^{r(\al,{\overline 0})}\right).
\end{equation*}
Similarly, we have
\begin{align*}
[S({\mathfrak L}_{\overline 1})]
&=\sum_{k=0}^{\infty}(-1)^k [S^k({\mathfrak L}_{\overline 1})] \\
&=\prod_{\al\in \G} \left(\sum_{m=0}^{\infty}(-1)^m [S^m({\mathfrak
L}_{(\al,\overline 1)})]e^{m(\al,\overline 1)}\right).
\end{align*}
By the definition of symmetric power operations, we have
\begin{equation*}
[S({\mathfrak L}_{\overline 1})]=\prod_{\al\in \G } \exp \left(
 \sum_{r=1}^{\infty}
\frac{{(-1)}^r}{r} \Psi^r([ {\mathfrak L}_{(\al,\overline 1)}])
e^{r(\al,\overline 1)} \right).
\end{equation*}
It follows that
\begin{align*}
&[C({\mathfrak L})]=\sum_{k=0}^{\infty}(-1)^k [C_k({\mathfrak L})]
=[\Lambda({\mathfrak L}_{\overline 0})]
\cdot [S({\mathfrak L}_{\overline 1})] \\
&=\prod_{\al \in \G} \exp
\left(-\sum_{r=1}^{\infty}\frac{1}{r}\Psi^r([{\mathfrak
L}_{(\al,\overline 0)}])
e^{r(\al,\overline 0)}\right)\\
& \times \prod_{\al \in \G}
\exp\left(\sum_{r=1}^{\infty}\frac{(-1)^r}{r} \Psi^r([{\mathfrak
L}_{(\al,\overline 1)}])e^{r(\al,\overline 1)} \right).
\end{align*}
Replacing $|{\mathfrak L}_{(\al,a)}|=(-1)^a[{\mathfrak L}_{(\al,a)}]$ and
$E^{\ag}={(-1)}^ae^{\ag}$, the above identity reduces to
\begin{equation*}
[C({\mathfrak L})]= \prod_{(\al,a)\in \G \times \mathbb Z_2}
\exp\left(- \sum_{r=1}^{\infty}\frac{1}{r}\Psi^r(|{\mathfrak
L}_{(\al,a)}|) E^{r(\al,a)}\right).
\end{equation*}
Hence, by the Euler-Poincar\'{e} principle, we obtain the {\it Lie
module denominator identity} for the graded Lie superalgebra
${\mathfrak L}=\bigoplus_{\ag \in \GG} {\mathfrak L}_{\ag}$.

\begin{prop}
For every $(\G \times \mathbb Z_2)$-graded Lie superalgebras
$\mathfrak L,$ we have
\begin{equation} \label {denominator}
\prod_{(\al,a)\in \G \times \mathbb Z_2} \exp\left(-
\sum_{r=1}^{\infty}\frac{1}{r}\Psi^r(|{\mathfrak L}_{(\al,a)}|)
E^{r(\al,a)}\right) =1-[H({\mathfrak L})].
\end{equation}
\end{prop}

\subsection{Main results on logarithmic functions}
To begin with, we introduce the $\C$-linear operators on
$\Gamma_K(G)[[\G \times \mathbb Z_2]]$ such as
\begin{equation}
\omega=\sum_{k\ge 1}\frac {\mu(k)}{k} \Theta_k \circ\Psi^k,\qquad
\eta=\sum_{k\ge 1}\frac 1k \Theta_k \circ\Psi^k,
\end{equation}
where $\mu$ is the M\"obius inverse function.
Recall that the maps $\Theta_k$ are defined
in Section \ref{Definition:Logarithm} and the maps
$$\Psi^k \,:\Gamma_K(G) [[\G\times \mathbb Z_2]]
\to \Gamma_K(G)[[\G\times \mathbb Z_2]]$$
are the induced algebra homomorphisms defined by
the action of Adams operations on coefficients.
Applying the M\"obius inversion formula one can easily show that
$\omega$ and $\eta$ are mutually inverses to each other.
Observe that
$$\eta([\mathfrak L ])=\displaystyle\sum_{(\gamma,a)
\in \G \times \mathbb Z_2}\eta(\gamma,a) E^{(\gamma,a)}$$
is the formal power series whose coefficients are given by
\begin{equation*}
\begin{split}
&\eta(\gamma,\overline 0)=\sum_{d|\gamma} \frac{1}{d} \Psi^{d}
\left(|\mathfrak L_{(\frac{\gamma}{d},\overline 0)}| \right )
+\sum_{d|\gamma \atop \text{ $d$: even }} \frac{1}{d}
\Psi^{d}\left(|\mathfrak L_{(\frac{\gamma}{d},\overline 1)}| \right ), \\
&\eta(\gamma,\overline 1)=\sum_{d|\gamma \atop \text{ $d$: odd }}
\frac{1}{d} \Psi^{d} \left(|\mathfrak
L_{(\frac{\gamma}{d},\overline 1)}| \right ).
\end{split}
\end{equation*}
Let
$$\mathcal D=\omega \circ\text{\rm Log}.$$
By definition $\omega$ is $\C$-linear and
$$\omega(R[[\G\times \mathbb Z_2]]^{(n)})
\subseteq R [[\G\times \mathbb Z_2]]^{(n)}.$$
Therefore, we conclude that $\mathcal D$ is logarithmic
by Proposition \ref {basis theorem}.

\begin{thm} \label {logarithm function}
$\mathcal D$ is the unique logarithmic function on
$\Gamma_K(G)[[\G \times \mathbb Z_2]]$ such that
\begin{equation}\label{property of log}
[\mathfrak  L]=\mathcal D([H(\mathfrak  L)]).
\end{equation}
for every $(\G \times \mathbb Z_2)$-graded Lie superalgebra
$\mathfrak L.$
Moreover $\omega \circ\text{\rm Log}=\text{\rm Log}\circ \omega$.
\end{thm}

\noindent{\bf Proof.}
If we take the logarithm on both sides of
the Lie module denominator identity \eqref {denominator}, then we
obtain the equality
\begin{equation}\label{take log}
\sum_{\ag}\eta(\al,a)E^{\ag}=\text{Log}([H(\mathfrak L)])
\end{equation}
where
$$\eta(\al,a)=\sum_{r|(\al,a)}\dfrac{1}{r}\Psi^r
\left(|\mathfrak L_{\frac{(\al,a)}{r}}|\right).$$
Equivalently,
\begin{equation*}
\eta([\mathfrak L])=\text{Log}([H(\mathfrak L)]).
\end{equation*}
Taking $\omega$ on both sides, we have
the desired result. On the other hand, since the operators
$\Theta_k$ and $\Psi^k$ commute with $\text{Log}$, so does
$\omega.$

In order to prove the uniqueness of $\mathcal D$ assume that
$\mathcal D'$ is another logarithmic function satisfying
Eq. \eqref{property of log}.
Let us choose a $\C$-basis of
$\G_K(G)[[\G \times \mathbb Z_2]]$, say
 $\{I_\ld E^{(\al,a)}\,:\,\ld \in \Lambda
,\,(\al,a)\in (\G \times \mathbb Z_2) \}$, such that
$I_\ld$ is an actual finite dimensional $KG$-module for every $\ld \in \Lambda$.
Given $\ld$ and $(\beta,b)$, let $V=\bigoplus_{(\al,a)\in \G \times \mathbb
Z_2} V_{(\alpha,a)}$, where
\begin{equation*}
V_{(\alpha,a)}=\begin{cases}
I_{\ld } & \text{ if }(\alpha,a)=(\beta,b),\\
0& \text{ otherwise.}
\end{cases}
\end{equation*}
Then $\mathcal D'$ must map $[V]$ to $[\mathfrak  L(V)]$
which coincides with $\mathcal D([V])$
by the property \eqref{property of log}.
In other words,
$$\mathcal D(I_\ld E^{(\beta,b)})=\mathcal D'(I_\ld E^{(\beta,b)})$$
for every $(\beta,b)\in \G \times \mathbb Z_2.$
It is obvious
$[\mathfrak L(V)]\in \G_K(G)[[\G \times \mathbb Z_2]]^{(ht (\beta))}.$
Such a logarithmic function is uniquely determined by Proposition
\ref {basis theorem} (a). So, $\mathcal D=\mathcal D'$.
\qed

Comparing the coefficients of both
sides of Eq. \eqref {take log}, we can derive that
$\eta(\al,a)$ equals to $W\ag$, which is defined
in \eqref{coef of logarithmic},
for all $\ag\in (\G\times \mathbb Z_2).$
It follows that
\begin{align*}
\mathcal D([\mathfrak L])
=\omega\circ \text{Log}([H(\mathfrak L)])
=\omega(\sum W\ag E^{\ag}).
\end{align*}
Immediately the above identity provide a close formula
for the homogeneous component $|{\mathfrak L}_{\ag}|$.

\begin{cor}\label {root space1 }
For every $\ag\in (\G\times \mathbb Z_2)$
\begin{equation}\label{close formula:Green algebra}
|{\mathfrak L}_{\ag}|=\sum_{ d > 0 \atop \ag=d(\tau,b)}
\frac{\mu(d)}{d}\Psi^d\left(W \left(\tau,b\right)\right),
\end{equation}
where
$$W(\alpha,a)=\sum_{s\in T(\alpha,a)}\dfrac{(|s|-1)!}{s!}
\prod {|H(\mathfrak L)_{(\beta_i,b_j)}|}^{s_{ij}}.$$
\end{cor}

\begin{cor}
\label {free lie superalgebra} Let $V=\bigoplus_{(\alpha,a)\in \G
\times \mathbb Z_2} V_{(\alpha,a)}$ be a $(\G \times \mathbb
Z_2)$-graded $K$-vector space with finite dimensional homogeneous
subspaces, and let $\mathfrak  L(V)=\bigoplus_{(\alpha,a)\in
\Gamma \times \mathbb Z_2} V_{(\alpha,a)}$ be the free Lie
superalgebra generated by $V.$ Suppose $G$ acts on $V$ preserving
$(\G \times \mathbb Z_2)$-gradation. Then,
$\mathcal D$ is the unique logarithmic function on $\Gamma(G)[[\G
\times \mathbb Z_2]]$ such that
\begin{equation}\label{free-case}
[\mathfrak  L(V)]=\mathcal D([V])
\end{equation}
for every $V.$
\end{cor}

\noindent{\bf Proof.}
By close inspection of the proof of
\cite[Corollary 3.2]{KK98}, we know that
$$H_k(\mathfrak  L(V))=\begin{cases}V & \text{ if }k=1,\\
0& \text{ otherwise }\end{cases}$$ for every field $K,$ char$K\ne
2.$
Hence, \eqref {property of log} is reduced to
\begin{equation*}
[\mathfrak  L(V)]=\mathcal D([ V]).
\end{equation*}
The uniqueness of $\mathcal D$ also follows
from Theorem \ref {logarithm function}.
\qed

For example, let us consider rank=1 case.
If we use the notation,
$[V]:=[V_{\overline 0}]-[V_{\overline 1}]$ and
$[{\mathfrak L(V)}_n]:=[{\mathfrak L (V)}_{(n,\overline
0)}]-[{\mathfrak L(V)}_{(n,\overline 1)}],$
then Eq. \eqref{free-case} implies that
\begin{equation*}
[\mathfrak  L(V)_{n}]=\dfrac 1n
\sum_{d|n}\mu(d)\Psi^d\left([V]^{\frac nd}\right).
\end{equation*}
Similarly, if we use the notation
$\{{\mathfrak L(V)}_n\}:=[{\mathfrak L (V)}_{(n,\overline
0)}]+[{\mathfrak L(V)}_{(n,\overline 1)}],$
then
Eq. \eqref {denominator} implies the formula
\begin{equation*}
\{\mathfrak  L(V)_{n}\}=\dfrac 1n
\sum_{d|n}\mu(d)\Psi^d\left([V_{\overline 0}]^{\frac nd}
-(-1)^d[V_{\overline 1}]^{\frac nd}\right).
\end{equation*}


From now on, we assume that $K$ is an arbitrary field with char$K\ne 2.$
Recall that we showed in
Theorem \ref {logarithm function} that if $\G_K(G)$ is
a special $\ld$-ring, then there exists a unique logarithmic function
$\mathcal D$ on $\Gamma_K(G)[[\G \times \mathbb Z_2]]$ satisfying
\begin{equation*}\label{property of log 1}
[\mathfrak  L]=\mathcal D([H(\mathfrak  L)])
\end{equation*}
for every $(\G \times \mathbb Z_2)$-graded Lie superalgebra
$\mathfrak L.$ But when $\G_K(G)$ is not a special $\ld$-ring, we
do not know whether this is true. However, in the case of free
Lie superalgebras, we can derive an analogue of Theorem \ref {logarithm
function}.
In showing this the subsequent lemma, which
can be obtained by applying the {\it Lazard elimination} theorem,
plays an essential role.

\begin{lem}{\rm (cf. \cite{Bour,B})}\label{lazard elimination theorem}
Let $G$ be a group and $K$ be any field with char$K\ne 2.$ For
$(\G \times \mathbb Z_2)$-graded $KG$-modules $U$ and $V,$ we have

{\rm (a)} $$\mathfrak L(U\oplus V)= \mathfrak L(U)\oplus \mathfrak
L(V\wr U),$$ where $V\wr U$ is the space spanned by all products
$[v,u_1,\cdots,u_m]$ with $m\ge 1, v\in V$ and $u_i\in U$ for all
$i$, with grading induced by that of $\mathfrak L(U\oplus V).$

{\rm (b)} $[V\wr U]=[V](1-[U])^{-1}.$
\end{lem}

\begin{thm} \label {green ring-logarithm function}
Let $G$ be a group and $K$ be any field with char$K\ne 2.$
Then, there exists a unique logarithmic function $\mathcal D$ on $
\G_K(G)[[\G \times \mathbb Z_2]]$ such that
\begin{equation}
[\mathfrak L(V)]=\mathcal D([V])
\end{equation}
for every $(\G \times \mathbb Z_2)$-graded $KG$-module $V.$
\end{thm}

\noindent{\bf Proof.}
Choose
$\{I_\ld E^{(\al,a)}\,:\,\ld \in \Lambda
,\,(\al,a)\in (\G \times \mathbb Z_2) \}$,
a $\C$-basis of
$\G_K(G)[[\G \times \mathbb Z_2]]$,
such that $I_\ld$ is an actual finite dimensional $KG$-module
for every $\ld \in \Lambda.$
Given $I_\ld E^{(\beta,b)}$, consider the graded $KG$-module
$V(\ld,(\beta,b))=\bigoplus_{(\al,a)\in \G \times \mathbb Z_2}
V_{(\alpha,a)}$, where $V_{(\beta,b)}=I_{\ld }$ and
$V_{(\al,a)}=0$ for all $(\al,a)\ne(\bt,b).$
Then, by Proposition
\ref {basis theorem} (b), there exists a unique logarithmic
function $\mathcal D$ such that
\begin{equation}\label{minimal case}
\mathcal D(I_\ld E^{(\beta,b)})
=[\mathfrak L(V(\ld,(\beta,b)))].
\end{equation}
It is obvious that
$[\mathfrak L(V(\ld,(\beta,b)))] \in \G_K(G)[[\G \times \mathbb
Z_2]]^{(ht (\beta))}.$
For any graded $KG$-module $V,$ write $[V]=\sum
a_{\ld,(\beta,b)}I_\ld E^{(\beta,b)}.$
We claim that $\mathcal
D([V]) =[\mathfrak L(V)].$
To show this observe that
\begin{equation}\label{second condition of log}
\mathcal D([U]+[V])=\mathcal D([U])+\mathcal D([V](1-[U])^{-1}).
\end{equation}
for $(\G \times \mathbb Z_2)$-graded $KG$-modules $U$ and $V$.
Indeed this property follows from the definition of
logarithmic function easily.
Combining \eqref {second condition of log}
with Eq. \eqref{minimal case} and Lemma \ref {lazard
elimination theorem}, we conclude that $\mathcal D([V])
=[\mathfrak L(V)].$
\qed

It is very worthwhile remarking that Theorem \ref {green ring-logarithm
function} guarantees the existence of such a
logarithmic function $\mathcal D,$ but provides no information on
its explicit form.
However, using Grothendieck algebra $\overline \Gamma_K(G)$
instead of Green algebra $\Gamma_K(G),$
we can write out such a logarithmic function explicitly.

From now on, let $G$ be a finite group and
$K$ be an arbitrary field.
Let $G_{p'}$ be the set of all elements of $G$ of order not
divisible by $p$.
And we let $\mathcal C$ be the $\mathbb C$-algebra
consisting of all class functions from $G_{p'}$ to $\mathbb C,$
that is, functions from $G_{p'}$ to $\mathbb C$ such that $f\in
\mathcal C,\,f(g)=f(g')$ whenever $g$ and $g'$ are conjugate in
$G$. Denote the algebraic closure of $K$ by $\hat K.$ We choose
and fix a primitive $e$-th root of unity $\xi$ in $\hat K$ and
$\omega \in \C,$ where $e$ denotes the least common multiple of
the orders of the elements of $G_{p'}.$ For a $KG$-module $V,$ we
define $\ch_{\overline V}$ be the function from $G_{p'}$ to $\C$
such that, for $a\in G_{p'},\, \ch_{\overline
V}(a)=\omega^{k_1}+\cdots+\omega^{k_r},$ where $\xi^{k_i}(\,1\le
i\le r) $ are eigenvalues of $a$ in its action on $\hat K \otimes
V$ (see \cite{B}).

\begin{lem}{\rm(\cite{Ben,B})}\hfill

{\rm (a)} There exists an injective $\C$-algebra homomorphism
$\tau: {\overline \G}_K(G) \to {\overline \G}_{\hat K}(G)$ such
that $\tau(\overline V)=\overline {\hat K \otimes V}.$

{\rm (b)} The $\C$-algebra homomorphism $\ch \, :\, {\overline
\G}_{\hat K}(G) \to \mathcal C$, defined by $ \ch(\overline I_\ld)
=\ch_{\overline I_\ld},$ is an isomorphism.
\end{lem}

Thus we may regard $\overline \G_{K}(G)$ as a subalgebra of
${\overline \G_{\hat K}(G)}.$
Then $\mathcal C$ becomes a special
$\Psi$-ring for the operations $\Psi^n(f)(g)=f(g^n)$ for all $g\in
G_{p'}$ and $n\ge 1.$ From Theorem \ref {psi=ld}, it
follows that $\mathcal C$ is a special $\ld$-ring.
Moreover, it is easy to
verify that
$$\ch_{\Psi^d(x)}=\Psi^d(\ch(x))$$
for $x \in \overline {\G}_K(G)$.
Also, we can easily show that $\overline
\G_{K}(G)$, when regarded as a subalgebra of $\mathcal C$, is
invariant under $\Psi^n$ for all $n.$

\begin{prop}\label{compact case}
Let $G$ be a finite group and $K$ be any field.
Then, $\overline
\G_K(G)$ is a special $\ld$-ring.
\end{prop}

\noindent{\bf Proof.}
Since $\overline \G_{K}(G)$, when regarded
as a subalgebra of $\mathcal C$, is invariant under $\Psi^n$ for
all $n,$ it becomes a special $\ld$-ring by Theorem \ref {psi=ld}.
\qed

\begin{rem}
In the case where $G$ is finite and $K$
has characteristic 0 or characteristic not dividing the order of
$|G|,$ it is well known that $\overline \G_{K}(G)$ may be identified
with ${\G}_{K}(G)$ {\rm (}see {\rm \cite {Ben,B}}{\rm )}. However, for $(G,K)$ such
that $G$ is finite and char$K|\,|G|,$ we have no idea whether
$\G_K(G)$ is a special $\ld$-ring or not.
\end{rem}

Consider a ($\G \times \mathbb Z_2$)-graded $K$-vector space
$V=\bigoplus_{(\alpha,a)\in (\G \times \mathbb
Z_2)}V_{(\alpha,a)}$ with $\dim V_{(\alpha,a)} < \infty$ for all
$(\alpha,a) \in \G \times\mathbb Z_2$. Suppose $G$ acts on $V$
preserving the $(\G\times \mathbb Z_2)$-gradation. We define
$$ \overline V:=
\sum_{(\alpha,a)} \overline {V_{(\alpha,a)}}e^{(\alpha,a)} \in
{\overline \G_{K}(G)}[[\G \times \mathbb Z_2]].$$ By setting
$|\overline {V_{(\alpha,a)}}|=(-1)^a \overline {V_{(\alpha,a)}},$
we have $ \overline V= \sum_{(\alpha,a)} |\overline
{V_{(\alpha,a)}}|E^{(\alpha,a)} .$

Let us consider a ($\G \times \mathbb Z_2$)-graded Lie
superalgebra $\mathfrak L=\bigoplus_{(\alpha,a)\in \G \times
\mathbb Z_2}\mathfrak L_{(\alpha,a)}$ with $\dim \mathfrak L_{(\alpha,a)}
< \infty$ for all $(\alpha,a) \in (\G \times\mathbb Z_2)$ over
$K$. Suppose $G$ acts on $\mathfrak L$ preserving the ($\G \times
\mathbb Z_2$)-gradation.
Then, we can derive
an identity analogous to the denominator identity Eq. \eqref {denominator}.
\begin{equation}\label{denomainator-Grothendieck}
\prod_{(\al,a)\in \G \times \mathbb Z_2} \exp\left(-
\sum_{r=1}^{\infty}\frac{1}{r}\Psi^r(|\overline {{\mathfrak
L}_{(\al,a)}}|) E^{r(\al,a)}\right) =1-\overline {H({\mathfrak L})}.
\end{equation}
We call this the {\it Grothendieck Lie module denominator identity}
of $\mathfrak L$.
Set
$$\overline {\mathcal D}=\omega \circ\text{\rm Log}.$$

\begin{thm}\label {char k-logarithm function}
$\overline {\mathcal D}$ is the unique logarithmic function on
$\overline \G_K(G)[[\G \times \mathbb Z_2]]$ such that
\begin{equation}\label{property of log II}
\overline {\mathfrak L}=\overline {\mathcal D}(\overline {H
(\mathfrak L)})
\end{equation}
for all every $(\G \times \mathbb Z_2)$-graded Lie superalgebra
$\mathfrak L.$
Moreover $\omega \circ\text{\rm Log}=\text{\rm Log}\circ \omega$.
\end{thm}

\noindent{\bf Proof.}
This can be exactly in the same way as in the proof of
Theorem \ref {logarithm function}.
\qed

We call $\overline {\mathfrak L}$ the {\it Grothendieck Lie module
function of $\mathfrak L.$} Comparing the coefficients of both
sides of Eq. \eqref {property of log II}
provides the following closed formula
analogous to Eq. \eqref{close formula:Green algebra}.
\begin{cor}\label {chark-root space1 }
For every $(\alpha,a) \in (\G \times\mathbb Z_2)$ we have
\begin{equation*}
\overline {|{\mathfrak L}_{\ag}|}=\sum_{ d > 0 \atop \ag=d(\tau,b)}
\dfrac {\mu(d)}{d}\Psi^d\left(W \left(\tau,b \right)\right),
\end{equation*}
where
$$W(\alpha,a)=\sum_{s\in T(\alpha,a)}\dfrac{(|s|-1)!}{s!}
\prod {\overline {|H(\mathfrak L)_{(\beta_i,b_j)}|}}^{s_{ij}}.$$
\end{cor}

\begin{cor} \label {char k-free logarithm function}
Let $V=\bigoplus_{(\alpha,a)\in \G\times \mathbb Z_2}
V_{(\alpha,a)}$ be a $(\G \times \mathbb Z_2)$-graded $K$-vector
space with finite dimensional homogeneous subspaces, and let
$\mathfrak  L(V)$
be the free Lie superalgebra generated by
$V.$ Suppose $G$ acts on $V$ preserving $(\G \times \mathbb
Z_2)$-gradation. Then, $\overline {\mathcal D}$ is the unique
logarithmic function on $\overline \G_K(G)[[\G \times \mathbb
Z_2]]$ such that
\begin{equation*}
\overline {{\mathfrak L}(V)}=\overline {\mathcal D}(\overline V)
\end{equation*}
for every $V.$
\end{cor}

\section{APPLICATIONS}\label{APPLICATIONS}

\subsection{New interpretation of the symmetric
power map, $\tilde s_t$, using plethysm}
\label{PLETHY}

Since Atiyah and Tall suggested
the "splitting principle" and ``verification principle",
the theory of special $\ld$-rings has been developed with
a close connection with the theory of symmetric functions
(\cite{AD}).
In particular, some identities and properties
of symmetric functions seem to be more natural
in our framework.
In this section, we reformulate plethystic equations
in the context of symmetric functions as relations
among the ring of Witt vectors, the necklace ring, and
the Grothendieck ring of formal power series.

\begin{exm}\label{example:section4}
Let $R$ be the ring of symmetric functions
in infinitely many variables $x_1,x_2,\cdots$.
We also let $e_n$ be the $n$-th elementary symmetric functions defined
using the product
$\prod_i(1+x_it),$
$h_n$ be the $n$-th complete symmetric functions defined using the product
$\prod_i\frac {1}{1-x_it},$ and
$p_n$ be the $n$-th power sum symmetric function.
Then, one can check that $R$ has a $\ld$-ring structure if we set
$\ld^n(e_1)=e_n$, or equivalently $\Psi^n(p_1)=p_n$ for all $n\ge 1.$
Observe that over the ring $R\otimes \mathbb Q$, the identity
\begin{align*}
\sum_{n=0}^{\infty}h_n t^n
=\prod_{n=1}^{\infty}\frac {1}{1-q_n t^n}
=\exp\left(\sum_{n=1}^{\infty}\frac {p_n}{n}t^n \right)
\end{align*}
can be rewritten as
\begin{align*}
&(q_1,q_2,\cdots)\stackrel{\Phi}{\mapsto}(p_1,p_2,\cdots),\\
&(h_1,h_2,\cdots)\stackrel{\text{int}^{-1}\circ
\frac {d}{dt}\log}{\mapsto}(p_1,p_2,\cdots).
\end{align*}
The symmetric function $q_n, (n\ge 1)$ enjoy the peculiar
property such that $(-q_n), n\ge 2,$ are Schur-positive
$($see {\rm \cite{ST}}$)$.
Similarly, if we define $(t_n)_n$ by
\begin{align*}
\sum_{n=0}^{\infty}e_n t^n
=\prod_{n=1}^{\infty}\frac {1}{1-t_n t^n}
=\exp\left(\sum_{n=1}^{\infty}\frac {(-1)^np_n}{n}t^n \right)
\end{align*}
$(-t_n), \,n\ge 2$ also turn out to be Schur-positive.
Indeed, this follows from the observation
\begin{equation*}
t_n=
\begin{cases}q_n & \text{ if }n \text{ is odd},\\
\text{sgn}_{S_n}\otimes q_n& \text{ if }n \text{ is even}.
\end{cases}
\end{equation*}
\end{exm}

Let the base ring be
$$\mathbb Q[[\Psi^m(d_n)\,:\,m,n\ge 1]].$$
Given the following commutative diagram
\begin{equation*}
\begin{CD}
(c_n)_n @>{\tau_R}>>
(b_n)_n@>{\tilde s_t}>> \displaystyle\sum_{n\ge 0}a_n t^n \\
@V{\Phi}VV @V{\tilde\varphi}VV  @V{\left(\frac {d}{dt}\right)\log }VV \\
(d_n)_n
@>id >> (d_n)_n
@>\text{\rm identification} >> \displaystyle\sum_{n\ge 0} d_{n+1} t^n,
\end{CD}
\end{equation*}
we let
$$A=\sum_{n\ge 1}a_n, \quad B=\sum_{n\ge 1}b_n,\quad
C=\sum_{n\ge 1}c_n.$$
Set the degree of $d_n$ to be $n$ for all $n\ge 1.$
Then, $A$ and $C$ can be viewed as the elements in
$$\mathbb Q[[d_j\,:\, j\ge 1]]$$
since $a_n$ and $c_n$ can be expressed as homogeneous polynomials
of degree $n$ with $\mathbb Q$-coefficients
in variables $d_j, 1\le j\le n$.
And
\begin{equation*}
B=\sum_{n=1}^{\infty}\dfrac 1n\sum_{d|n}\mu(d)\Psi^d(d_{\frac nd})
\end{equation*}
since
\begin{equation*}
b_n=\dfrac 1n\sum_{d|n}\mu(d)\Psi^{d}(d_{\frac nd}).
\end{equation*}
In order to generalize the notion of {\it plethysm}
for elements
$$\xi(d_j\,;\, j\ge 1) \in \mathbb Q[[d_j\,:\, j\ge 1]]$$
and
$$\zeta(\Psi^{i}(d_j)\,;\, i,j\ge 1)\in \mathbb Q[[\Psi^m(d_n)\,:\,m,n\ge 1]]$$
let us define $\xi \circledcirc \zeta$ by
\begin{equation}\label{def of plethysm}
\xi\circledcirc \zeta:=\xi(\zeta_1,\zeta_2,\cdots),
\end{equation}
where
$$\zeta_k=\zeta(\Psi^{ki}(d_j)\,;\, i,j\ge 1).$$
With this notation
we can obtain another characterization of $\tilde s_t$.
\begin{thm}\label{plethysm}
{\rm (a)} With the above notation, we have
$$A \circledcirc B=A.$$
{\rm (b)}
Set the degree of $\Psi^n(d_j)$ to be $nj$ for all $n,j\ge 1$.
Let $Y=\sum_{n= 1}^{\infty}y_n,$ where $y_n$ is a homogeneous polynomial in
$\Psi^d(d_{\frac nd})$ with $d, n\ge 1, \text{ and }d\,|\,n$.
Then, $A\circledcirc Y=A$ if and only if $Y=B.$
Equivalently, $A\circledcirc Y=A$ if and only if
$\tilde s_t(y_1,y_2,\cdots)=1+\sum_n a_nt^n.$
\end{thm}
\noindent{\bf Proof.}
(a) By the definition \eqref{def of plethysm} we obtain
\begin{align*}
A\circledcirc B &=\exp\left(\sum_{n=1}^{\infty}\frac {d_{n}}{n} \right)\circledcirc B\\
&=\exp\left(\sum_{n=1}^{\infty}\frac {B_n}{n} \right).
\end{align*}
Here,
\begin{equation*}
B_n=\sum_{i=1}^{\infty}\dfrac 1i\sum_{d|i}\mu(d)\Psi^{nd}(d_{\frac id}).
\end{equation*}
On the other hand,
\begin{align*}
\sum_{n=1}^{\infty}\dfrac{B_n}{n}&=\sum_{n=1}^{\infty}\sum_{m\ge 1 \atop d|m}
\dfrac 1n \dfrac 1m  \mu(d)\Psi^{dn}(d_{\frac md})\\
&=\sum_{s,t}\Psi^{s}(d_t)\sum_{d|s}\dfrac {1}{st}\mu(d)\qquad
(\text{letting }dn=s,\,\frac md=t) \\
&=\sum_{t=1}^{\infty}\dfrac {d_t}{t}.
\end{align*}
Since $$A=\exp\left(\sum_{n=1}^{\infty}\frac {d_n}{n} \right),$$
the desired result follows.

(b) From the proof of (a) it follows that
$$\sum_{d|n}\dfrac 1d\Psi^d(y_{\frac nd})=\dfrac {d_n}{n}.$$
Applying M\"obius inversion formula implies that $y_n=b_n$ for all $n\ge 1$.
\qed

Consider the case $d_n=p_n$ for all $n\ge 1$.
In this case, $A\circledcirc B$ coincides with the usual plethysm of $A$ and $B$,
denoted by $A\circ B$.
Let $H=\sum_{n\ge 0}h_n$ and $L=\sum_{n\ge 1}l_n$,
where
\begin{align*}
l_n=\dfrac 1n\sum_{d|n}\mu(d){p_d}^{\frac nd}.
\end{align*}
Viewing $H$ and $L$ as functions in variables $p_n$'s, $n\ge 1$,
Theorem \ref{plethysm} (a) implies that
$$H \circ L=\dfrac {1}{1-p_1}.$$
Indeed, Joyal proved this identity using PBW theorem.
On the contrary, Reutenauer showed it using the plethysm (refer to \cite {R}).

Let
$$L(k)=\sum_{n=1}^{\infty}\dfrac 1n\sum_{d|n}c(k,d)\Psi^d(d_{\frac nd}),$$
where $c(k,d)$ is the sum of the $k$-th powers of the primitive $d$-roots of unity.
Then, by the same argument as in \cite{ST}, we have
$$A\circledcirc L(k)=\prod_{l|k}\dfrac {1}{d_l}.$$

Imitating the proof of Theorem \ref{plethysm}
we can provide some interesting plethystic equations.
Let
\begin{equation*}
\tilde l_n=\dfrac 1n\sum_{d|n}(-1)^{\frac nd}\mu(d){p_d}^{\frac nd}.
\end{equation*}
Regarding $H$, $L$, $E=\sum_{n\ge 0}e_n$,
and $\tilde L=\sum_{n\ge 1}\tilde l_n$ as functions
in $p_n$'s, $n\ge 1$,
we can provide the following plethystic equations
in the same way as Reutenauer did.
\begin{cor}\hfill

{\rm (a) } $H \circ \tilde L=\dfrac {1}{1+p_1}$\,\,\,\,\,\,\,
{\rm (b) } $E \circ L=\dfrac {1-p_1}{1-p_2}$   \,\,\,\,\,\,\,
{\rm (c) } $E \circ \tilde L=\dfrac {1+p_1}{1+p_2}$
\end{cor}
\noindent{\bf Proof.}
Since the identities (a) throughout (c) can be proven in the same way,
we will prove only (c).
\begin{align*}
E\circ L &=\exp\left(\sum_{n=1}^{\infty}(-1)^n\frac {p_{n}}{n} \right)\circ \tilde L\\
&=\exp\left(\sum_{n=1}^{\infty}(-1)^n\frac {{\tilde L}_n}{n} \right).
\end{align*}
Thus, by definition of $\tilde l_n$
we have
\begin{align*}
\sum_{n=1}^{\infty}(-1)^n \dfrac{\tilde L_n}{n}&=\sum_{n=1}^{\infty}\sum_{m\ge 1 \atop d|m}
(-1)^n (-1)^{\frac md}  \dfrac 1n \dfrac 1m  \mu(d)p_{dn}^{\frac md}\\
&=\sum_{s,t}(-1)^t p_s^t\sum_{d|s}\dfrac {1}{st}(-1)^{\frac sd}\mu(d)
\qquad (\text{letting }dn=s,\,\frac md=t) \\
&=-\sum_{t=1}^{\infty}\dfrac {p_1^t}{t}+ \sum_{t=1}^{\infty}\dfrac {p_2^t}{t}\\
&=\log \left(\dfrac {1+p_1}{1+p_2}\right).
\end{align*}
The third equality follows from
$$\sum_{d|n}\mu(n)(-1)^{\frac nd}=
\begin{cases}-1 &\quad\text{\rm if }n=1 \\
2 & \quad\text{\rm if }n=2\\
0 & \quad \text{\rm otherwise. }
\end{cases}$$
Therefore, we have the desired result.
\qed

\subsection{Generators of supersymmetric functions}
\label{supersymmetric*}

In this section, we are going to show that
the diagram \eqref{big diagram}, being applied to supersymmetric polynomials,
provides several generating sets of the set of supersymmetric polynomials.
We recall the definition of supersymmetric polynomials briefly (see \cite{Ste}).
Let $K$ be a field of characteristic $0$,
and let $x_1,\cdots,x_a,y_1,\cdots,y_b,t$ be independent indeterminates.
A polynomial $p$ in
$$K[\mathbb X, \mathbb Y]:=K[x_1,\cdots,x_a,y_1,\cdots ,y_b]$$
is said to be supersymmetric if

(1) $p$ is invariant under permutations of $x_1,\cdots,x_a,$\\
(2) $p$ is invariant under permutations of $y_1,\cdots,y_b,$\\
(3) when the substitution $x_1=t,\,y_1=t$ is made in $t$, the resulting
polynomial is independent of $t$.

Let $T(a,b)$ denote the set of supersymmetric polynomials in
$K[\mathbb X, \mathbb Y].$
In \cite{Ste} Stembridge provided two generating sets of $T(a,b)$
such as
$$\{\sigma^{(n)}_{a,b}(x;y): n\ge 1\}$$
and
$$\{\tau^{(n)}_{a,b}(x;y): n\ge 1 \}.$$
Here, $\sigma^{(n)}_{a,b}(x;y):=(x_1^n\cdots +x_a^n)-(y_1^n\cdots +y_b^n)$
and $\tau^{(n)}_{a,b}(x;y)$ is defined by the equation
$$\sum_{n=0}^{\infty}\tau_{a,b}^{(n)}(x;y)t^n
:=\prod_{i,j} \frac{1-x_it}{1-y_jt}.$$
Now, let us consider the element $r-s \in R$ where $r$ is $a$-dimensional
and $s$ $b$-dimensional.
After decomposing $r,s$ into the sum of $1$-dimensional elements
we can write  $r=x_1+\cdots+x_a$ and $s=y_1+\cdots+y_b$.
Considering the following diagram
\begin{equation*}
\begin{CD}
(q_{a,b}^{(n)})_n @>{\tau_R}>>
(r,0,0,\cdots)@>{\tilde s_t}>> \displaystyle\sum_{n\ge 0}h_{a,b}^{(n)} t^n \\
@V{\Phi}VV @V{\tilde \varphi}VV  @V{\left(\frac {d}{dt}\right)\log }VV \\
(\Psi^n(r))_n
@>id >> (\Psi^n(r))_n
@>\text{\rm identification} >> \displaystyle\sum_{n\ge 0} \Psi^{n+1}(r) t^n,
\end{CD}
\end{equation*}
we obtain supersymmetric polynomials
$q_{a,b}^{(n)}(x;y)$ and $h_{a,b}^{(n)}(x;y)$ satisfying
\begin{equation}\label{generating set0}
\begin{aligned}
\sum_{n=0}^{\infty} h_{a,b}^{(n)}(x;y) t^n
&=\prod_{n=1}^{\infty}\frac {1}{1-q_{a,b}^{(n)}(x;y) t^n}\\
&=\exp\left(\sum_{n=1}^{\infty}\frac {\sigma_{a,b}^{(n)}(x;y)}{n}t^n \right).
\end{aligned}
\end{equation}
Note that $h^{(n)}_{a,b}(x;y)$ is defined by the equation
$$\sum_{n=0}^{\infty}h_{a,b}^{(n)}(x;y)t^n :=\prod_{i,j} \frac{1-y_jt}{1-x_it}.$$
Similarly, we can consider supersymmetric polynomials
$t_{a,b}^{(n)}(x;y)$ and $e_{a,b}^{(n)}(x;y)$ such that
\begin{equation}\label{generating set1}
\begin{aligned}
\sum_{n=0}^{\infty} e_{a,b}^{(n)}(x;y) t^n
&=\prod_{n=1}^{\infty}\frac {1}{1-t_{a,b}^{(n)}(x;y) t^n}\\
&=\exp\left(\sum_{n=1}^{\infty}\frac {(-1)^n\sigma_{a,b}^{(n)}(x;y)}{n}t^n \right),
\end{aligned}
\end{equation}
where $e^{(n)}_{a,b}(x;y)$ is defined by the equation
$$\sum_{n=0}^{\infty}e_{a,b}^{(n)}(x;y)t^n :=\prod_{i,j} \frac{1+x_it}{1+y_jt}.$$
Finally, we let
\begin{equation}\label{supersymmetric function l}
l_{a,b}^{(n)}(x;y)=\dfrac 1n\sum_{d|n}\mu(d){\sigma_{a,b}^{(d)}(x;y)}^{\frac nd}.
\end{equation}
In view of Eq. \eqref{generating set0} through Eq. \eqref{supersymmetric function l}
we can obtain several generating sets.

\begin{prop}
The sets of supersymmetric polynomials
$\{h^{(n)}_{a,b}(x;y):n\ge 1\}$, $\{e^{(n)}_{a,b}(x;y): n\ge 1\}$,
$\{q^{(n)}_{a,b}(x;y):n\ge 1\}$, $\{t^{(n)}_{a,b}(x;y):n\ge 1\}$,
and $\{l^{(n)}_{a,b}(x;y):n\ge 1$
generate the algebra $T(a,b)$ respectively.
\end{prop}

\subsection
{Recursive formulas for $|\mathfrak L_{(\al,a)}|$ and $|\overline
{\fL_{(\al,a)}}|$} \label{RECURSIVE FORMULA}

In Corollary \ref {root space1 } and Corollary \ref
{chark-root space1 } we provided closed formulas which
enable us to compute the homogeneous components
$|\mathfrak L_{(\al,a)}|$ and $|\overline
{\fL_{(\al,a)}}|$ for a graded Lie superalgebra $\mathfrak L$ if
we know information on the homology $H(\mathfrak L).$
On the contrary, we focus on recursive formulas
associated with (Grothendieck)Lie module denominator identities
in this section.
For this end, we introduce some formal
$R$-linear operators on $R [[\G\times \mathbb Z_2]],$ where $R$ is
a commutative algebra with unity over $\C.$

Let $\widehat{\Gamma}_{\C} = \C \otimes_{\Z} \widehat{\Gamma}$ be
the complexification of $\widehat {\Gamma}$. Choose a
non-degenerate symmetric bilinear form $(\ | \ )$ on
$\widehat{\Gamma}_{\C}$ and fix a pair of dual bases
$\{u_i\}\text{ and } \{u^i\}$ of $\widehat{\Gamma}_{\C}$. We
define a partial ordering $\ge$ on $\Gamma_{\C}$ by $\lambda \ge
\mu $ if and only if $\lambda-\mu \in \Gamma$ or $\lambda =\mu$.
We will denote by $\lambda > \mu$ if $\lambda \ge \mu$ and
$\lambda \neq \mu$. In particular, if $\zeta(\la,a)\in \C,$ it
means that $\zeta(\la,a)\cdot 1.$

\begin{df}{\rm (\cite{KKO})}\hfill

{\rm (a)} The {\it partial differential operators } $\partial_i$
and $\partial^i$ are defined by
$$\partial_i(E^{(\lambda,a)})=(\lambda|u_i)E^{(\lambda,a)}, \quad
\partial^i(E^{(\lambda,a)})=(\lambda|u^i)E^{(\lambda,a)}.$$

{\rm (b)} For an element $\rho\in\widehat{\Gamma}_{\C}$, we define
the {\it $\rho$-directional derivative} $D_{\rho}$ by
\begin{equation*}
D_{\rho}(E^{(\lambda,a)})
=\sum(\rho|u_{i})\partial^{i}(E^{(\lambda,a)})=(\rho|\lambda)E^{(\lambda,a)}.
\end{equation*}

{\rm (c)} The {\it Laplacian} $\Delta$ is defined to be
\begin{equation*}
\Delta(E^{(\lambda,a)})
=\sum\partial^{i}\partial_{i}(E^{(\lambda,a)})
=(\lambda|\lambda)E^{(\lambda,a)}.
\end{equation*}
\end{df}

Recall the definition of $\eta$ given in Section 3.1.
Define
$$\eta^{*}([\mathfrak L])= \displaystyle \sum_{(\gamma,a) \in \G \times \mathbb Z_2} \eta^{*}(\gamma, a)
E^{(\gamma,a)}$$ to be the formal power series whose coefficients
are given by
\begin{equation*}
\eta^{*}(\gamma,a)=(\gamma| \gamma)\eta(\gamma,a) -
\sum_{(\gamma,a) =(\gamma',a')+( \gamma'',a'')} (\gamma'|\gamma'') \eta(\gamma',a') \eta(\gamma'',a'').\\
\end{equation*}
Similarly, we write
$$\eta(\overline {\mathfrak L} )=\displaystyle\sum_{(\gamma,a) \in \G \times \mathbb Z_2} \overline
{\eta} (\gamma,a) E^{(\gamma,a)}$$ and
$$\eta^{*} (\overline {\mathfrak L} )=\displaystyle\sum_{(\gamma,a) \in \G \times \mathbb Z_2} \overline
{\eta^{*}} (\gamma,a) E^{(\gamma,a)}.$$

With this notation,
\begin{prop}\label {recursive formula}
For $(G,K)$ such that char$K\ne 2$ and $\G_K(G)$ is a special
$\ld$-ring, we have \hfill

{\rm (a)} $D_{\rho}(1-[H(\mathfrak L)])=-D_{\rho}(\eta([ \mathfrak
L]))(1-[H(\mathfrak  L)]).$

{\rm (b)} $\Delta(1-[H(\mathfrak L)])=-\eta^*([ \mathfrak L
])(1-[H(\mathfrak  L)]).$
\end{prop}

\noindent{\bf Proof.}
By applying the following formal
differential identities to the Lie module denominator identity
\eqref {denominator}
\begin{equation}\label{diff-operator}
\begin{aligned}
&D_{\rho}(\text{log}D)=\frac {D_{\rho}(D)}{D},\\
& \frac{\Delta (D)}{D}=\sum_i \partial_i \left(\frac{\partial^i
D}{D}\right) +\sum_i \left(\frac{\partial_i D}{D}\right)
\left(\frac{\partial^i D}{D}\right)
\end{aligned}
\end{equation}
and then comparing the coefficients of both sides, we can complete
the proof.
\qed

Recall that under the condition that $G$ is a finite group and $K$
a field with char$K\ne 2,$ we obtained the Grothendieck Lie module
denominator identity \eqref {denomainator-Grothendieck} for every
$(\G\times \mathbb Z_2)$-graded Lie superalgebra $\mathfrak L$. By
applying \eqref {diff-operator} to Eq. \eqref{denomainator-Grothendieck},
we can derive the following relations.

\begin{prop}\label {recursive formula 1}
Let $G$ be a finite group and $K$ be any field with char$K\ne 2.$
Then, we have \hfill

{\rm (a) $ D_{\rho}(1- \overline {H(\mathfrak
L)})=-D_{\rho}(\eta(\overline {\mathfrak L}))(1-\overline {
H(\mathfrak L)}).$

(b) $\Delta(1-\overline {H(\mathfrak  L)})=-\eta^*(\overline {
\mathfrak L})(1-\overline {H(\mathfrak  L)}).$ }

\end{prop}

\begin{rem}
If $\mathfrak  g$ is a Borcherds
superalgebra, the operator $(\Delta +2D_{\rho})$ plays an
essential role. In fact, in this case, we have
$$(\Delta +2D_{\rho})(1-[H(\mathfrak  g_-)])=0,$$
where $\rho$ is a Weyl vector. Exploiting this fact, we can obtain
Peterson's and Freudenthal's formulas for the $[{\mathfrak
g}_{\alpha}].$
\end{rem}

Comparing the coefficient of $E^{(\al,a)}$ in Proposition \ref
{recursive formula} (a) and (b), we have

\begin{cor}\label{recur1}
\begin{align*}
& {\rm (a)}\quad(\rho|\alpha) \eta(\alpha,a) - \sum_{\beta <
\alpha} (\rho|\beta)
\eta(\beta,b) |H(\mathfrak L)_{(\alpha-\beta,a-b)}| =(\rho|\alpha)| H(\mathfrak L)_{(\alpha,a)}|. \\
&{\rm (b)}\quad   \eta^*(\alpha,a)- \sum_{\beta < \alpha
}\eta^*(\beta,b) |H(\mathfrak L)_{(\alpha-\beta,a-b)}|
=(\alpha|\alpha) |H(\mathfrak L)_{(\alpha,a)}|.
\end{align*}
\end{cor}

Similarly, it follows from Proposition \ref {recursive formula 1}
(a) and (b) that
\begin{cor}
\begin{align*}
& {\rm (a)}\quad (\rho|\alpha) \overline {\eta}(\alpha,a) -
\sum_{\beta < \alpha} (\rho|\beta)\overline{ \eta}(\beta,b)
|\overline {H(\mathfrak L)_{(\alpha-\beta,a-b)}}| =(\rho|\alpha)
|\overline { H(\mathfrak L)_{(\alpha,a)}}|. \\
&{\rm (b)}\quad  \overline{\eta^*}(\alpha,a)- \sum_{\beta < \alpha
}\overline{\eta^*}(\beta,b) |\overline {H(\mathfrak
L)_{(\alpha-\beta,a-b)}}| =(\alpha|\alpha) |\overline {H(\mathfrak
L)_{(\alpha,a)}}|.
\end{align*}
\end{cor}

For example, let $M$ be the Monster simple group and
$\fL=\bigoplus_{(m,n)}\fL_{(m,n)}$
be the Monster Lie algebra (see \cite{B92}).
Recall that $\fL_{(m,n)}\simeq V_{mn}$ for $(m,n)\ne 0$, where
$V^\natural=\bigoplus_{n\ge -1}V_{n}$ is the Moonshine Module
constructed by Frenkel et al.
Here, we use $\G$-grading, not
the $\G\times \mathbb Z_2$-grading.
Applying Corollary \ref{recur1} (a),
we can recover the identity
$$[\fL_{(m,n)}]=-\sum_{k|(m,n)\atop k>1}
\dfrac 1k \Psi^k\left([\fL_{(\frac mk,\frac nk)}]\right)+
\sum_{1\le k<m\atop 1\le l<n}\dfrac km c(k,l)[H(\fL)_{(m-k,
n-l)}]+[H(\fL)_{(m, n)}],$$ where $c(m,n)=\displaystyle \sum_{k\ge
1\atop k|(m,n)}\dfrac 1k \Psi^k\left([\fL_{(\frac mk,\frac
nk)}]\right).$

\subsection{Replicable functions from the viewpoint of
logarithmic functions}
\label{REPLICABLE}

The concept of replicable functions was first introduced
by Norton as a generalization of the {\it replication formulae}.
In their famous {\it Moonshine conjecture} Conway and Norton
suggested replication formulae as an important
family of character identities that are satisfied
by the Thompson series of the {\it Monster simple group}
(see \cite{N}).
Later, Borcherds has proven this conjecture completely
by showing that the Thompson series are indeed replicable functions.

Let $F(q) =q^{-1}+\sum_{n\ge 1}f(n)q^n$ be a normalized
$q$-series. By setting $q = e^{2 \pi i z}$ with $\text{Im} z >0$,
we often write $F= F(z)$ so that the notation can be consistent
with the Fourier expansion of modular functions. Note that for
each $m \ge 1$, there exists a unique polynomial $X_m(t)\in
\mathbb C[t]$ such that
\begin{equation*}
X_m(F)\equiv \frac 1m q^{-m} \quad \text{mod} \ q\mathbb C[[q]].
\end{equation*}
We write
\begin{equation*} \label{eq:H_{m,n}}
X_m(F)=\frac{1}{m} q^{-m}+\sum_{n \ge 1}H_{m,n}q^n.
\end{equation*}
In \cite{N}, Norton has shown that the coefficients $H_{m,n}$ satisfy
the identity
\begin{equation}\label{eq:Norton1}
\sum_{m,n\ge1}H_{m,n}p^m
q^n=-\log\left(1-pq\sum_{i=1}^{\infty}f(i)\frac{p^i-q^i}{p-q}\right).
\end{equation}
Indeed, one can show that Eq. \eqref{eq:Norton1}
is equivalent to the product identity
\begin{equation}\label {eq:Norton2}
p^{-1}\prod_{m=1}^{\infty} \exp \left(-X_m(F(q))p^m \right)
=F(p)-F(q).
\end{equation}

We recall the definition of {\it replicable functions} (\cite
{ACMS,B92,N}).

\begin{df} {\rm
A normalized $q$-series $F(q)=q^{-1}+\sum_{n\ge 1}f(n)q^n$ is said
to be {\it replicable} if $H_{a,b}=H_{c,d}$ whenever $ab=cd$ and
$(a,b)=(c,d)$.}
\end{df}

The replicable functions can be characterized as follows.
\begin{prop}{\rm (\cite {ACMS,N})}
A normalized $q$-series $F(q) =q^{-1}+\sum_{n\ge 1}f(n)q^n$ is
replicable if and only if for all $m>0$ and $a|m,$ there exist
normalized $q$-series $F^{(a)}(q)=q^{-1} +
\sum_{n=1}^{\infty}f^{(a)}(n)q^n$ such that
\begin{equation*}\label{eq:replicable1}
F^{(1)}=F, \quad X_m(F)=\frac 1m\sum_{ad=m \atop 0 \le b <d}
F^{(a)}\left(\frac{az+b}{d}\right).
\end{equation*}
where $q=e^{2\pi i z},\, Im\,z >0. $
\end{prop}

The normalized $q$-series $F^{(a)}$ is called the {\it $a$-th
replicate} of $F$. If $F^{(a)}$ are also replicable for all $a\ge
1$, then $F$ is said to be {\it completely replicable}.

Let $\Gamma=\Z_{>0}\times \Z_{>0}$ and consider the formal power
series ring
$$\C[[\Gamma]]=\{\sum_{m,n \ge 0,
\atop (m,n)\ne (0,0)}a(m,n)p^mq^n|\, a(m,n)\in \C\}.$$
Set
$$T(m,n)=\{s={(s_{ij})}_{i,j\ge 1}|\,s_{ij}
\in \mathbb Z_{\ge 0},\sum s_{ij}(i,j)=(m,n)\},$$
and define
$$W(m,n):=\sum_{s\in T(m,n)}\dfrac{(|s|-1)!}{s!}\prod _{i,j}f(i+j-1)^{s_{ij}}.$$
Then, we have
\begin{equation}\label{eq:Norton 3}
\text{Log}\left(\sum_{i,j\ge 1 }f(i+j-1)p^i
q^i\right)=\sum_{m,n\ge1}W(m,n)p^m q^n.
\end{equation}

\begin{prop}A normalized $q$-series $F$ replicable if and only
if $W(a,b)=W(c,d)$ whenever $ab=cd$ and $(a,b)=(c,d).$
\end{prop}

\noindent{\bf Proof.}
From Eq. \eqref {eq:Norton2} one can derive that
\begin{equation}\label{h(m,n)=w(m,n)}
\prod_{m,n=1}^{\infty} \exp \left(-H_{m,n}p^mq^n \right)
=1-\sum_{m,n=1}^{\infty}f(m+n-1)p^m q^n.
\end{equation}
If we take the logarithm on both sides of Eq. \eqref{h(m,n)=w(m,n)}, then
it follows from \eqref {eq:Norton 3} that $H_{m,n}=W_{m,n}.$ This
completes the proof.
\qed

From now on, we discuss how to obtain completely-replicable
functions in a unified way.

First, we consider monstrous functions
appearing in Moonshine conjecture.
The Monster Lie algebra $\mathfrak L$ is a $II_{1,1}$-graded
representation of the Monster simple group $M$ acting by
automorphisms of $\mathfrak L$ for $(m,n) \ne (0,0)$ as
$M$-modules. In particular,
$$\text{tr}(g|\mathfrak L_{(m,n)})=\text{tr}(g|V_{mn}):=c_g (mn)\text{ for }
 g\in M,\quad(m,n)\ne(0,0),$$
where $V$ is the {\it Moonshine module} $V^\natural
=\bigoplus_{n\ge -1}V_n$ constructed by Frenkel et al. and $c_g
(n)$ is the coefficient of $q^n$ of the elliptic modular function
$J(q)=j(q)-744.$ From the denominator identity of $\mathfrak L$ in
\cite {B92}, we obtain
\begin{equation*}
\prod_{m,n \ge 1}\exp \left(-\sum_{d|(m,n)}\dfrac 1d
\Psi^d([V_{\frac{mn}{d^2}})])p^mq^n\right)=1-\sum_{i,j}[V_{i+j-1}]p^iq^j.
\end{equation*}

Let
$$\G(M)[[\Gamma]]=\{\sum_{m,n \ge 0,\atop (m,n)\ne (0,0)}a(m,n)p^mq^n|\, a(m,n)\in
\G(M) \}.$$ Define $\mathcal D: \G(M)[[\G]] \to \G(G) [[\Gamma]]$
by $\mathcal D=\omega \circ \text{\rm Log}$, where $\omega=\sum_{k\ge
1}\dfrac {\mu(k)}{k} \Theta_k \circ\Psi^k.$ Then, we have

\begin{prop}\label{characterization OF MONSTER}
$\mathcal D$ is a logarithmic function on $\G(M)[[\G]]$ such that
\begin{equation}\label{property of MONSTER}
\sum_{m,n\ge 1}[V_{mn}]p^mq^n=\mathcal D\left(\sum_{m,n\ge
1}[V_{m+n-1}]p^m q^n \right).
\end{equation}
\end{prop}

Taking $\ch_g$ on both sides the above identity \eqref{property of
MONSTER}, we have
\begin{equation}\label{property of MONSTER1}
\sum_{m,n\ge 1}\text{\rm tr}(g|V_{mn})p^mq^n=\mathcal
D\left(\sum_{m,n\ge 1}\text{\rm tr}(g|V_{m+n-1})p^m q^n \right),
\end{equation}
which says that Thompson series $T_g$ are replicable functions.

As for replicable functions, in particular non-monstrous
functions, the above method is no more effective. In these cases
we need a different approach.

For every positive integer $r,$ consider the formal power series
in $q$
$$h^{(r)}=q^{-1}+\sum_{m=1}^{\infty}x_m^{(r)}q^m$$
whose coefficients are the indeterminates $x_m^{(r)}$. Let
$$\Xi =\C[\cdots,x_m^{(r)},\cdots \,|\,m,r\ge 1].$$
For a fixed $r$ and an arbitrary $m \ge 1$
consider the family of equations
\begin{equation}\label{formally replicable}
X_m(h^{(r)}(q))-\frac 1m \sum_{ad=m \atop 0 \le b <d}
h^{(ra)}\left(\exp(2\pi i \frac bd)q^{\frac ad} \right).
\end{equation}
Expand \eqref {formally replicable} in a $q$-series and then
consider the coefficient of $q^n$ for every $n\ge 1.$ Let $I^{(r)}$ be
the ideal in $\Xi$ generated by them.  Let I be the ideal in $\Xi$
generated by $\bigcup_{r=1}^{\infty}I^{(r)}.$

Consider the semigroup $\Delta$ of $GL^+(2,\,\mathbb Q)$ given by
$$\Delta=\left \{ \begin{pmatrix} a&b\\0&d \end{pmatrix} |\,a,b,d \in \mathbb Z,\, a>0
\right \}.$$ For $A=\begin{pmatrix}a&b\\0&d\end{pmatrix}\in
\Delta,$ set
$$e||_A=e,\quad x_m^{(n)}||_A=x_m^{(rn)},\quad q||_A=\exp(2\pi i \frac bd)q^{\frac ad}$$
for all $e\in \C, \, m,m\ge 1.$ We extend the mapping $||_A$ to
the $\C$-algebra homomorphism from $\Xi [[q]]\to \Xi [[q]].$ Since
the ideal $I$ is stable under $||_A,$ it induces a homomorphism
from $\Xi/I$ to $\Xi/I.$

\begin{df}
Let $n$ be the positive integer. We define $\Psi^n:\Xi/I \to
\Xi/I$ to be
$$\Psi^n(x)=x||_A$$
for all $x\in \Xi/I$ where $A=\begin{pmatrix} n&0\\0&n
\end{pmatrix} .$
\end{df}

\begin{prop}{\rm
$\Xi/I$ is a special $\ld$-ring. }
\end{prop}

\noindent{\bf Proof.}
It is easy to verify that
$\Xi/I$ is a special $\Psi$-ring. By Theorem
\ref {psi=ld}, we conclude that $\Xi/I$ is a special $\ld$-ring.
\qed

Define $\mathcal D: \Xi /I [[\Gamma]] \to \Xi /I [[\Gamma]]$ by
$\mathcal D=\omega \circ \text{\rm Log},$ where $p=e^{(1,0)},$
$q=e^{(0,1)}$ and $\omega=\sum_{k\ge 1}\dfrac {\mu(k)}{k} \Theta_k
\circ\Psi^k.$ Then, we have

\begin{prop}\label{characterization OF REP I}
$\mathcal D$ is a logarithmic function on $\Xi/I[[\G ]]$ such that
\begin{equation}\label{property of REPLICABLE}
\sum_{m,n\ge 1}x^{(r)}_{mn}p^mq^n=\mathcal D\left(\sum_{m,n\ge
1}x_{m+n-1}^{(r)}p^mq^n\right)
\end{equation}
for every positive integer $r.$
\end{prop}

\noindent{\bf Proof.}
By the identity \eqref {eq:Norton2} we have
\begin{equation*}\label {eq:Norton3}
p^{-1}\prod_{m=1}^{\infty} \exp \left(-X_m(h^{(r)}(q))p^m \right)
=h^{(r)}(p)-h^{(r)}(q)
\end{equation*}
for all $r\ge 1.$ Substituting $$\dfrac 1m \displaystyle\sum_{ad=m
\atop 0 \le b <d} h^{(ra)}\left(\exp(2\pi i \frac bd)q^{\frac ad}
\right)$$
 for
$X_m(h^{(r)}(q))$ gives rise to the following product identity
\begin{equation*}
\prod_{m,n=1}^{\infty} \exp \left(-\sum_{k=1}^{\infty} \frac{1}{k}
\Psi^{k}(x^{(r)}_{mn}) p^{km} q^{kn}\right)
=1-\sum_{m,n=1}^{\infty}x^{(r)}_{m+n-1}p^m q^n.
\end{equation*}
Taking the logarithmic function $\mathcal D$ on both sides, we get
Eq. \eqref {property of REPLICABLE}.
\qed

Let $F(q)=q^{-1} + \sum_{n=1}^{\infty}f(n)q^n$ be a completely
replicable function and $F^{(a)}(q)=q^{-1} +
\sum_{n=1}^{\infty}f^{(a)}(n)q^n$ be its $a$-th replicates for all
$a\ge 1.$ Then, we obtain a $\C$-algebra homomorphism, $\psi_F
:\Xi/I \to \C$, such that $\psi_F(x_n^{(a)})=f^{(a)}(n).$
Conversely, if we have a $\C$-algebra homomorphism,
$\psi :\Xi/I \to \C,$ we get a completely-replicable function
$F_\psi^{(1)}$ and its replicates
by setting the q-series $F_\psi^{(r)}(q):=q^{-1}+\sum_{n\ge
1}\psi(x_n^{(r)})q^n$. Define $\Phi$ by the function
from the set completely replicable functions to the set of
$\C$-algebra homomorphisms from $\Xi/I$ to $\C$ sending $F$ to
$\psi_F,$ and $\Upsilon$ by the function from the set of
$\C$-algebra homomorphisms from $\Xi/I$ to $\C$ to the set
completely replicable functions sending $\psi$ to $F_\psi.$ Then,
we have

\begin{prop}\label{characterization OF REP II}
There is a natural one-to-one correspondence between the set of
completely replicable functions and the set of $\C$-algebra
homomorphisms from $\Xi/I$ to $\C.$
\end{prop}

\noindent{\bf Proof. }
It suffices to show that $\Phi\circ \Upsilon=id$ and
$\Upsilon\circ \Phi=id,$ equivalently
$$F_{\psi_F}=F \text{ and } \psi_{F_\psi}=F,$$
which follows from the definition of $F_\psi$ and $\psi_F$
immediately.
\qed

In the above proposition, if we replace $I$ by $I^{(1)},$ then we
have
\begin{prop}\label{characterization OF REP III}
There is a natural one-to-one correspondence between the set of
replicable functions and the set of $\C$-algebra homomorphisms
from $\Xi/I^{(1)}$ to $\C.$
\end{prop}

\noindent{\bf Acknowledgements}

The author would like to express his
sincere gratitude to Professor Kenji Iohara  and Soon-Yi Kang
for their interest in this work and many valuable advices.

\end{document}